\documentclass[12pt]{puthesis_onarxiv}

\usepackage{amssymb}
\usepackage[dvips]{graphics}

\newcommand{\map}{\mbox{$f:(T^2,\tau)\rightarrow S^3$ }}
\newcommand{\proj}{\mbox{$\pi:X\rightarrow \mathbb{C}P^1$}, \mbox{$(x,y)\mapsto x$ }}
\newcommand{\ra}{\rightarrow}

\newcommand{\Wpm}{\Omega_{\pm}}
\newcommand{\Wp}{\Omega_+}
\newcommand{\Wm}{\Omega_-}
\newcommand{\Cpm}{C_{\pm}}

\newcommand{\cp}{\mbox{$\mathbb{C}P^1$}}
\newcommand{\Wpmh}{\mbox{$\widehat{\Omega}_{\pm}$ }}
\newcommand{\Wph}{\widehat\Omega_+}

\newcommand{\pR}{\frac{\partial}{\partial R}}
\newcommand{\pmu}{\frac{\partial}{\partial\mu}}
\newcommand{\pnu}{\frac{\partial}{\partial\nu}}
\newcommand{\ppnu}{\frac{\partial^2}{\partial\nu^2}}
\newcommand{\pli}{\frac{\partial}{\partial\lambda_i}}

\newcommand{\m}{\mu}
\newcommand{\n}{\nu}
\newcommand{\im}{\sqrt{-1}}
\newcommand{\be}{\begin{eqnarray}}
\newcommand{\ee}{\end{eqnarray}}
\newcommand{\bes}{\begin{eqnarray*}}
\newcommand{\ees}{\end{eqnarray*}}
\newcommand{\p}{p^\prime}
\newcommand{\po}{\p_0}
\newcommand{\ka}{\kappa}
\newcommand{\w}{\omega}
\newcommand{\Oone}{O\left(\frac{1}{\mu}\right)}
\newcommand{\R}{\mathbb{R}}
\newcommand{\ds}{\displaystyle}
\newcommand{\z}{\zeta}

\newcommand{\lieG}{{\mathfrak g}}
\newcommand{\x}{x}
\newcommand{\ppr}{p^\prime}

\title{On the Existence of Minimal Tori in $S^3$ of Arbitrary Spectral Genus.}
\author{Emma Elizabeth Carberry}
\abstract{Harmonic mappings are a generalisation of geodesics, and are defined as the solutions to a natural variational problem. Interest in them began in
1873 with Plateau's problem of finding surfaces of minimal area bounded by
given closed space curves. The field has been studied by
mathematicians and physicists ever since, and is now both broad and extremely active. In this dissertation I consider harmonic maps which can be studied using integrable
systems, and thus by algebro-geometric means. In particular I focus upon a simple case of both geometric and physical interest, namely harmonic maps $f$ from a 2-torus
(with conformal structure $\tau$) to the 3-sphere. In \cite{hitchin:90} Hitchin showed that (except in the case of a conformal map to a totally geodesic $S^2\subseteq S^3$) the data $(f,\tau)$ is in one-to-one correspondence with certain algebro-geometric data. This data consists of a hyperelliptic curve $X$ (called the spectral curve) together with a projection map $\pi:X\ra \mathbb{C}P^1$, a pair of holomorphic functions on $X-\pi^{-1}\{0,\infty\}$, and a line bundle on $X$, all satisfying certain conditions. He proved (case-by-case) that for $g\leq 3$, there are curves of genus g that support the required data, and hence describe harmonic maps $f:(T^2, \tau)\rightarrow S^3$. Once such a curve is found, the line bundle
may be chosen from a real $(g-2)$-dimensional family, and each choice yields a new harmonic map. Of especial interest are conformal harmonic maps as their images are minimal surfaces. I show that for each $g\geq0$, there are countably many conformal harmonic maps $f:(T^2, \tau)\rightarrow S^3$ whose spectral curves have genus~$g$. All of these harmonic tori have rectangular conformal type, and those with $g>2$ come in real $(g-2)$-dimensional families. %
}
\acknowledgements{I have been fortunate to receive support from a number of people during my years in Princeton. It is a pleasure to acknowledge them here. My advisor, Phillip A. Griffiths, has been a generous source of support and mathematical wisdom. I have benefited greatly from his skilled guidance and superb insight. The many discussions I had with Ivan Sterling during a visit to the University of Toledo helped to place my project in a larger perspective, which was important. Chuu-Lian Terng and Ian McIntosh also imparted to me excellent advice and suggestions. I met regularly with my fellow student Julianna Tymoczko for rigorous and detailed discussions, which proved most helpful.

My parents, sister Josie, grandmother and Uncle Roger did not allow a detail as small as the Pacific Ocean to obstruct their love and support, and by these my life has always been enriched. Last but not least, I wish to thank Dargan Frierson and my friends at Princeton, whose presence made my time here  ever so much more enjoyable.}
\dedication{\begin{center}{\it To Mum and Dad}\end{center}}

\newtheorem{lemma}{Lemma}[chapter]
\newtheorem{prop}{Proposition}[chapter]
\newtheorem{thm}{Theorem}[chapter]

\begin{document}
\chapter{Introduction}
A map $f:M\rightarrow N$ between Riemannian manifolds is {\em harmonic} if it is critical for the energy functional \mbox{$E(f)=\int_M \|df\|_p^2 dM$}, where $\|~\|_p$ is the norm on \mbox{$T_p^*M \otimes T_{f(p)}N$} induced from the Riemannian metrics on $M$ and $N$.  Examples of harmonic maps include geodesics, harmonic functions and holomorphic and antiholomorphic maps between K\"{a}hler manifolds. If $f$ is an isometric immersion, then it is critical for the energy functional if and only if it is critical for the area functional; hence the relevance of harmonic maps to Plateau's problem.  When $M$ is a surface, the harmonic map equations are invariant under conformal changes in the domain metric. Thus a conformal map of $M$ is harmonic precisely when its image is a minimal surface. The conformal invariance also means that one can consider harmonic maps of a Riemann surface.

Harmonic maps from a Riemann surface to a
compact Lie group or symmetric space are of particular interest. One reason for this is their
relationship with the important Yang-Mills equations of mathematical
physics. The harmonic map equations are then locally the self-dual
Yang-Mills equations on $\mathbb{R}^4$ with signature $(+,+,-,-)$, invariant under translation in the last two variables. Physicists study harmonic maps of $\mathbb R^2$ and
$\mathbb R^{1,1}$ in order
to gain insight into the Yang-Mills equations. Harmonic maps of surfaces also arise in the study of surfaces of
geometric interest. For example, the theorem of Ruh and Vilms~\cite{RV:70} asserts
that a surface has constant mean curvature precisely when its Gauss map is
harmonic. Notice that the Gauss map is then a harmonic map to $S^2$, and hence to $S^3$. Similar characterizations exist for both Willmore surfaces and surfaces of constant negative Gauss curvature~\cite{Bryant:84,MS:93}.



The last twenty-five years have seen an explosion of interest in this
area. 
A major theme of this research has been the  ``classification'' of harmonic maps. A series of papers (e.g.~\cite{Calabi:67,EW:83,Wolfson:88,Wood:88,Uhlenbeck:89,BR:90}) gave descriptions of harmonic maps from $S^2$ to various symmetric spaces in terms of an algebraic curve in an auxiliary complex manifold. Hitchin~\cite{hitchin:90} showed that harmonic maps from a 2-torus (with some complex structure $\tau$) to the 3-sphere (with standard metric) also enjoy an algebro-geometric description. Underlying this is the idea \cite{Pohlmeyer:76,ZM:78,ZS:79,Uhlenbeck:89} that by the insertion of a parameter into the harmonic map
equations, they can be reformulated as a family of equations of a
particularly pleasant form, namely as equations of Lax type. These equations linearise on the Jacobian of a hyperelliptic curve, called the {\em spectral curve}~\cite{Griffiths:85, Burstall:92}. Thus to a harmonic map $f:(T^2,\tau)\rightarrow S^3$ there corresponds a hyperelliptic curve. Hitchin also proved that if one begins with a hyper\-elliptic curve with certain additional data, one can construct a torus $(T^2,\tau)$  and a harmonic map $f:(T^2,\tau)\rightarrow S^3$. Hence in order to find harmonic tori in the 3-sphere, one seeks hyperelliptic curves supporting the requisite additional data. Not surprisingly, this is not a trivial problem. Hitchin proved (case-by-case) that for each $g\leq 3$, there is a torus $(T^2,\tau)$ and a harmonic map $f:(T^2,\tau)\rightarrow S^3$ whose spectral curve has genus $g$.  I show that for each $g\geq 0$, there are {\it conformal} harmonic maps of tori to $S^3$ whose  spectral curves have genus $g$. The importance of conformality is that the conformal harmonic maps are precisely those whose images are minimal surfaces. The precise algebro-geometric statement is given below. 

Given a curve $X$: \mbox{$y^2=x\prod_{i=1}^{g}{(x-\alpha_i)(x-\bar{\alpha_i}^{-1})}$},~\mbox{$|{\alpha_i}|\neq 0,1$}, let
\begin{itemize}
\item[--] $\pi$ be the projection to the $x$-plane,
\item[--] $\sigma$ be the hyperelliptic involution \mbox{$(x,y)\mapsto (x,-y)$},
\item[--] $\rho$ be the antiholomorphic involution,
\mbox{$(x,y)\mapsto(\bar{x}^{-1},(\prod_{i=1}^{g}{\lambda_i \bar{\lambda_i}^{-1}})^{\frac{1}{2}}\frac{\bar{y}}{{\bar{x}}^{g+1}})$}, where we choose the square root so that $\rho$ fixes the points with $|x|=1$,
\item[--] $\gamma_1$ be a curve in $X$ joining the two points in $\pi^{-1}(1)$, and
\item[--] $\gamma_{-1}$ a  curve in $X$ joining the two points in $\pi^{-1}(-1)$.
\end{itemize}
Then for each genus $g > 0$, there are countably many curves $X$ as above possessing meromorphic differentials
$\Theta$ and $\Psi$ that satisfy the following conditions:
\begin{enumerate}
\item $\Theta$ and $\Psi$ have double poles at $\pi^{-1}\{0,\infty\}$, and are holomorphic elsewhere. The principal parts of $\Theta$ and $\Psi$ are linearly independent over $\mathbb R$. (The ratio of these gives the conformal type $\tau$ of the torus.)
\item  $\Theta$ and $\Psi$ satisfy the symmetry conditions $\sigma^*\Theta=-\Theta$, $\sigma^*\Psi=-\Psi$, and the reality conditions $\rho^*\Theta=\bar{\Theta}$, $\rho^*\Psi=\bar{\Psi}$.
\item The integrals of $\Theta$ and $\Psi$ over $\gamma_1$, $\gamma_{-1}$, and over a basis for the homology of $X$, are integers.
\end{enumerate}

The last condition is by far the most difficult to satisfy; it demands that transcendental objects (the integrals) be integers. This places a severe restriction on the hyperelliptic curve.

My approach to this problem follows the work of Ercolani, Kn\"{o}rrer and Trubowitz~\cite{EKT:93}, in which they prove that for each even $g\geq 2$, there is a constant mean curvature torus in $\mathbb{R}^3$ whose spectral curve has genus $g$.

\chapter{Algebro-Geometric Description}
In his paper \cite{hitchin:90}, Hitchin showed that specifying a 
conformal structure $\tau$ on $T^2$, and a harmonic map 
$f:(T^2,\tau)\ra SU(2)$ (other than a branched conformal map to a totally 
geodesic $S^2\subseteq S^3$) is equivalent to specifying certain 
algebro-geometric data, described in theorems \ref{thm:general} and 
\ref{thm:ftns} below. He in fact 
gave an algebro-geometric description of harmonic sections of an 
$SU(2)$ principal bundle over $(T^2,\tau)$, and then specified the 
conditions the data must satisfy in order to correspond to the 
special case of a harmonic mapping. The purpose of this section is to 
give an expository account of Hitchin's work, outlining the ideas, 
but omitting most proofs.

Given  Riemannian manifolds $(M,g)$ and $(N,h)$, the {\it energy} 
$E(f)$ of $f:M\ra N$ is 
$$E(f)=\frac 12 \int_M \| df \|^2_p\, dM,$$
where $\|\cdot\|$ is the metric on $T^*_p M\otimes T_{f(p)}M$ induced 
from $g$ and $h$. We say that $f$ is {\it harmonic} if whenever 
$\{f_t\}$ is a one-parameter family of smooth maps with $f_0=f$,
$$\left.\frac{d}{dt}\right|_{t=0}E(f_t) =0.$$ 
If $M$ is a surface, then the harmonicity of $f$ depends only on the 
conformal class of the metric $g$, and so we can speak of harmonic 
maps of Riemann surfaces. If in addition $f:M\ra N$ is conformal, 
then it is harmonic if and only if $f(M)$ is a minimal surface. 

Let $M$ be a compact Riemann surface, $G$ a compact Lie group with 
bi-invariant metric,
and $f:M\rightarrow G$ a harmonic map. We denote the Lie algebra 
$T_eG$ of $G$ by $\lieG$. On $G$ we define the Maurer-Cartan form 
$\omega$ to be the unique $\lieG$-valued 1-form on $G$ that is 
left-invariant and acts as the identity on $T_e G$. $\omega$ 
satisfies the Maurer-Cartan equation 
$$d\omega + \frac {1}{2}[\omega,\omega]=0,$$ and for a linear group 
$\omega$ is given by $g^{-1}dg$. Given any smooth map $f$ from $M$ to 
$G$, we may pull back the Maurer-Cartan form to obtain a 
$\lieG$-valued 1-form $\phi$ on $M$ satisfying
\begin{equation}
d\phi+\frac 12 [\phi,\phi]=0.
\label{eq:MC}
\end{equation} 
Conversely, if $U$ is a 
simply-connected and connected open subset of $M$, then given a 
$\lieG$-valued 1-form $\phi$ on $U$ satisfying the Maurer-Cartan 
equation, we may integrate it to obtain a smooth map $f:U:\ra G$ such 
that $f^*(\omega)=\phi$, where $f$ is defined only up to left 
translation. Since the metric on $G$ is left invariant, one expects
that there is a condition characterising those $\phi$ that 
correspond to harmonic maps, and we calculate it below. 

For a one-parameter family  $\{f_t\}$,
\begin{eqnarray*}
\left.\frac{d}{dt}\right|_{t=0} E(f_t)&=&
\frac 12 \int_M \left.\frac{d}{dt}\right|_{t=0}\|df_t\|_p ^2\,dM\\
&=&\frac 12 \int_M \left.\frac{d}{dt}\right|_{t=0}\|df_t f_t^{-1}\|_p 
^2\, dM\\
&=&\frac 12 \int_M \left.\frac{d}{dt}\right|_{t=0}\langle 
df_tf_t^{-1} ,*(df_t f_t^{-1}) \rangle_p\,dM.
\end{eqnarray*}
Writing $f_t=e^{th}f$,
\begin{eqnarray*}
\left.\frac{d}{dt} E(f_t)\right|_{t=0}&=&
\int_M  \langle f dh f^{-1},*f^{-1}(df)\rangle_p\\
&=&
\int_M  \langle dh ,*(f^{-1}df)\rangle_p\\
&=&
\int_M  d\langle h ,*(f^{-1}df)\rangle_p -\langle h 
,d*(f^{-1}df)\rangle_p,
\end{eqnarray*}
so we see that $f$ is harmonic if and only if 
\begin{equation}
d*\phi=0.\label{eq:intparts}
\end{equation}

Equations~(\ref{eq:MC}) and (\ref{eq:intparts}) may be rewritten in 
terms of natural connections on $f^*TG$. $TG$ has connections $d_L$ 
and $d_R$ corresponding to its trivialisations by left and right 
translation, respectively. These are related by
$$d_R = d_L + ad_\omega\mbox{, (where $ad_\omega(a)=[\omega,a]$)}$$ 
and the Levi-Civita connection $d_A$ is the average of the two:
$$d_A = \frac 12 (d_L + d_R).$$ We will use the same notation for the 
pull-back of these connections to $f^*TG$, and we may rewrite 
(\ref{eq:MC}) as 
\begin{equation}
d_A(\phi)=0.\label{eq:dA}
\end{equation}

Denoting the $(1,0)$-component of $\frac 12 \phi$ by $\Phi$, we have
$$\frac 12 \phi = \Phi - \Phi^*,$$
and since
\begin{equation}
*(\Phi - \Phi^*)=\im(\Phi + \Phi^*),\label{eq:decomp}\end{equation} 
we have $$[\phi, *\phi]=0,$$ so that (\ref{eq:intparts}) is 
equivalent to 
\begin{equation}
d_A(*\phi)=0.
\label{eq:dAstar}
\end{equation}
Using (\ref{eq:decomp}) again, equations~(\ref{eq:dA}) and 
(\ref{eq:dAstar}) are 
\begin{eqnarray*}
\bar\partial_A\Phi - \partial_A\Phi^*&=&0\\
\bar\partial_A\Phi + \partial_A\Phi^*&=&0
\end{eqnarray*}
so give
$$\bar\partial_{A}€\Phi = 0.$$
We have also that $d_L=d_A - \frac{1}{2}\phi$ is flat, and hence
\begin{eqnarray*}
0&=&d_L^2=(d_A-\frac 12 \phi)^2=d_A^2+(\frac 12\phi)^2\mbox{, (from 
(\ref{eq:dA}))}\\
&=&F_A+(\Phi-\Phi^*)^2,
\end{eqnarray*}
so 
$$F_A=[\Phi,\Phi^*].$$

These equations may be considered as taking place in a trivial 
principal $G$-bundle over $M$, as follows. The trivial principal 
$G$-bundle 
$$\begin{array}{ccc}
G&\times&G\\
&\downarrow&\\
&G&
\end{array}$$
with projection $(g_1,g_2)\mapsto g_1$ and action $(g_1,g_2)\cdot 
h=(g_1,g_2 h)$ has the obvious trivial connection, which we shall 
denote by 
$\nabla_L$, along with the trivial connection $\nabla_R = 
g^{-1}\nabla_L g$. More explicitly, writing $\omega_L$ and $\omega_R$ 
for the connection $1$-forms of $\nabla_L$ and $\nabla_R$, then \begin{eqnarray*}
\omega_L &=& g_2^{-1}dg,\\
\omega_R &=& (g_1 g_2)^{-1}d(g_1 g_2)\\
&=&\omega_L + Ad_{g_2^{-1}}\omega.
\end{eqnarray*} 
We set $$\nabla_A=\frac{1}{2}(\nabla_L + \nabla_R).$$
Pulling back under the map $f$, we obtain a trivial principal bundle 
$P$, with connections $\nabla_L,\,\nabla_R,\,\nabla_A$ and 
$\phi=Ad_{g_2}(\omega_R-\omega_L)$. Note that $Ad(P)=f^*(TG)$. Again 
we may write $\frac{1}{2}\phi=\Phi^-\Phi^*$, and can consider $\Phi$ 
to be a section of $Ad(P)\otimes K$, where $K$ is the canonical 
bundle of holomorphic $1$-forms on $M$. 


Hitchin seeks solutions to the pair of equations \\
\parbox{13.8cm}{\begin{eqnarray*}
\bar\partial_A\Phi &=&0\\
F_A&=&[\Phi,\Phi^*]
\end{eqnarray*}}\hfill
\parbox{1cm}{\begin{eqnarray}\label{eq:}\end{eqnarray}}\\
for a connection $A$ in a principal $G$-bundle $P$ over $M$, and 
section $\Phi$ of $Ad(P)\otimes K$, called a {\it Higgs field}.  If 
these data come from a 
harmonic map $f:M\ra G$, then one additionally has that the 
connections
\begin{eqnarray*}
\nabla_L&=&\nabla_A-\Phi+\Phi^*\\
\nabla_R&=&\nabla_A+\Phi-\Phi^*
\end{eqnarray*}
are trivial. This is also a sufficient condition for the data 
$(A,\Phi)$ to come from a harmonic map as if $\nabla_L$, $\nabla_R$ 
as defined 
above are trivial then the difference between these trivialisations 
is a map $f:M\ra G$, and equations~(\ref{eq:}) tell us that $f$ is 
harmonic. It is however fruitful to first study general solutions to 
this pair of equations, and then focus on those corresponding to 
harmonic maps. The equations~(\ref{eq:}) are locally 
the self-dual Yang-Mills equations on $\mathbb{R}^4$ with signature 
$(+,+,-,-)$, invariant under translation in the last two variables.
The equations~(\ref{eq:}) can be made more conducive to analysis by 
the introduction of a parameter $x\in\mathbb{C}^*$, termed the 
{\it spectral parameter}. An easy check reveals that 
equations~(\ref{eq:}) are equivalent to the statement that for every 
$x\in\mathbb{C}^*$, the connections
\begin{equation}
d_x = d_A + x^{-1}\Phi - x\Phi^* \label{eq:family}
\end{equation}
are flat.
  
We shall henceforth narrow our focus to the case where $M$ is a torus 
with some conformal structure $\tau$ and $G=SU(2)$.  Then 
(\ref{eq:family}) defines a family of flat $SL(2,\mathbb{C})$ 
connections. $SU(2)\cong S^3$ is a natural and simple case of 
interest. Ultimately we will assign an algebraic curve to each 
solution of  equations~(\ref{eq:}); the fact that $SU(2)$ is a rank 
one symmetric space enables us to describe the solution using a 
single curve, and the fact that 
the elements of $SL(2,\mathbb{C})$ are $2\times 2$ matrices gives 
that this curve is hyperelliptic. We refer the reader to 
\cite{BFPP:93} for a study of harmonic maps of complex \mbox{$n$-tori} into 
symmetric spaces, and note in particular their result that any 
non-conformal harmonic map of a (real) 2-torus into a rank one 
symmetric space is of finite type. The restriction to tori is however 
more essential; an integrable systems approach to harmonic maps of 
higher genus Riemann surfaces has so far proven elusive. In what 
follows we will make use of the fact that the fundamental group 
$\pi_1(T^2)$ is abelian. 

Several geometrically interesting properties of a harmonic map $\map$ 
can be described in terms of the connection $A$ and Higgs field 
$\Phi$. In particular, $f$ is branched conformal if and only if 
$\det\Phi=0$ (with branch points the zeros of $\Phi$) and $f$ maps to a totally geodesic $S^{2}\subseteq 
S^{3}$ if and only if $A$ is reducible to a $U(1)$ connection. This 
latter condition is also equivalent to the existence of a gauge 
transformation $g$ leaving $A$ invariant and satisfying $g^{2}=-1$, 
$g^{-1}\Phi g=-\Phi$.

We have a holomorphic family of flat $SL(2,\mathbb{C})$ connections 
on $T^2$, and so we study the holonomy of these connections. We 
shall consider them as connections in a vector bundle $V$ over 
$(T^{2},\tau)$ 
with the $SU(2)$ structure exhibited by a quaternionic structure
$$j:V\ra V,\; j^{2}=-1,$$
and a symplectic form 
$$\omega:V\times V\ra \mathbb{C}.$$
Then $\Phi$ is a holomorphic section of End\,$V\otimes K$ with trace zero.
 
In the case where $\det\Phi =0$, $\Phi\not\equiv 0$, there are additional holomorphic 
invariants that we may associate to the solution. We then have 
$\Phi^{2} = 0$, and may define a holomorphic line bundle $L\subseteq 
V$ by $L\subseteq\mbox{ker\,}\Phi$. Thus $\Phi:L^{*}\cong V/L \ra 
L\otimes K$, so we may consider $\Phi$ as a holomorphic section $u$ 
of $L^{2}\otimes K$. It is natural to ask whether $d_{A}$ preserves the 
sub-bundle $L$; the obstruction being the holomorphic 
section $v$ of $L^{-2}\otimes K$ defined by 
$$v s^{2} = \omega(\partial_{A} s,s),$$ where $s$ is any local 
holomorphic section of $L$. If $v\not\equiv 0$, then the 
extension of $L$ defining $V$ is non-trivial. The quadratic 
differential $uv\in H^{0}(T^{2},K^{2})\cong\mathbb{C}$ \label{defn:uv} 
is then a constant multiple of the second 
fundamental form of $f(T^{2})\subseteq S^{3}$.
Assume additionally that $A$ is irreducible. Then $v\in H^0(T^2,L^{2}\otimes K)$ is not identically zero, so deg\,$L^2\leq 0$. If deg\,$L^2<0$, then $u$ vanishes identically, and hence so does $\Phi$. But from $F_A=[\Phi, \Phi^*]$ we have then that $A$ is flat, so is given by a representation of $\pi_1(T^2)$ in $SU(2)$. Since $\pi_1(T^2)$ is abelian, we see that $A$ reduces to a $U(1)$ connection. Thus deg\,$L^2=0$, and $u$ is nowhere vanishing. Hence a branched conformal harmonic map $f:(T^2,\tau)\ra S^3$ whose image does not lie in a totally geodesic $S^2$ is in fact a conformal immersion.  

Take $p\in T^{2}$ and generators $a$ and $b$  for 
$\pi_1(T^2,p)=\mathbb{Z}\otimes\mathbb{Z}$ that are conjugate to 
$[0,1]$ 
and $[0,\tau]$ respectively. For each $\x\in\mathbb{C}^*$, 
denote by $H(\x,p)$ and $K(\x,p)$ the holonomy of 
$d_\x$ around $a$ and $b$ respectively, and write 
$h(\x)={\rm tr}H(\x,p)$,  $k(\x)={\rm tr}K(\x,p)$. Note that 
the conjugacy classes of $H(\x,p)$ and $K(\x,p)$ are 
independent of the choice of base point $p$. Denoting by $\mu(\x)$ 
and $\mu^{-1}(\x)$ the eigenvalues of $H(\x,p)$, 
$$\mu^2(\x)-h(\x)\mu(\x)+1=0,$$
and so
$$\mu(\x)=\frac{1}{2}\left(h(\x)\pm\sqrt{h(\x)^2-4}\right).$$

This defines a 2-sheeted branched cover of $\mathbb{C}^{*}$, which we 
shall 
compactify to a cover of $\mathbb{C}P^{1}$. The resulting 
curve is algebraic:
\begin{prop}
For each $\x\in\mathbb{C}^{*}$, let $H(\x)$ be the holonomy 
of $d_{\x}$ around $[0,1]$, and $h(\x)={\rm tr} H(\x)$.
The function $h(\x)^{2}-4$ has finitely many odd order zeros 
in $\mathbb{C}^{*}$.
\end{prop}

Unsurprisingly, the proof of this result employs the fact that an 
elliptic operator on a compact domain has but a finite-dimensional 
kernel. It is worth remarking that it also utilises the full 
two-dimensional compactness of the torus, rather than the mere 
one-dimensional compactness of the closed curve $[0,1]$. 

If the holonomy is  identically trivial, $H(\x)\equiv 1,\, 
K(\x) \equiv 1$ (possibily after tensoring with a flat 
$\mathbb{Z}_{2}$-bundle), then this construction will not yield a 
two-sheeted covering, and cannot be described using the methods of 
\cite{hitchin:90}. However this occurs if and only if the data 
$(A,\Phi)$ correspond to a conformal map into a totally geodesic 
$S^{2}\subseteq S^{3}$, and such mappings can be described using 
divisors on $(T^{2},\tau)$. We shall assume henceforth that we 
are not in 
this trivial case.

Importantly, the branched cover of $\mathbb{C}^{*}$ defined by the 
holonomy $H(\x)$ is the same as that defined by $K(\x)$, 
this being is a consequence of the fact that the fundamental group 
of the torus $T^2$ is abelian. The next step is to study the behaviour 
of the holonomy of $d_{\x}$ as $\x\to 0$ and $\infty$, 
and thus determine the branching behaviour of the compactified curve. 
 
Now $A$ is an $SU(2)$ connection, and thus 
$j^{-1}d_{A}j=d_{A}$.
Since $j^{-1}\Phi j=-\Phi^{*}$, this gives 
$$j^{-1}d_{\x}j=d_{\bar\x^{-1}},$$
so $$\left(H(\x)^{-1}\right)^{*}=j^{-1}H(\x)j = H(\bar\x^{-1}),$$
and hence $$h(\bar\x^{-1})=\overline{h(\x)}.$$
Thus the behaviour of the holonomy as $\x\ra\infty$ is 
determined by that 
as $\x\ra 0$. Moreover, one can show that the eigenspaces of the 
holonomy 
matrices are determined simply by the holomorphic structure of $V$, 
and then use the regularity of this structure to study the holonomy 
in the limit $\x\ra 0$.
\begin{prop}
Let $(A,\Phi)$ be a solution of (\ref{eq:}) with 
$\det(\Phi)=-\eta^{2}dz^{2}\neq 0$. 
For each $\x\in\mathbb{C}^{*}$, we use $\mu(\x)$ and 
$\nu(\x)$ to denote 
the eigenvalues of the holonomy of $d_{\x}$ around $[0,1]$ 
and $[0,\tau]$ respectively. There is a punctured neighbourhood of 
$0\in\mathbb{C}$ such that
\begin{eqnarray*}
\pm\log\mu(\x) &=& \eta\x^{-1}+a+\x b(\x)\\
\pm\log\nu(\x) &=& \eta\tau\x^{-1}+\tilde{a}+\x 
\tilde{b}(\x) \end{eqnarray*} where $b(\x)$ and $\tilde b(\x)$ are even 
holomorphic functions.
\end{prop}

\begin{prop}
Let $(A,\Phi)$ be a solution of (\ref{eq:}) with 
$\det(\Phi)=0$ and $A$ irreducible. 
There is a punctured neighbourhood of 
$0\in\mathbb{C}$ such that
\begin{eqnarray*}
\pm \log{\mu(\x)} &=& \kappa\x^{-\frac{1}{2}}+\im 
k\pi+\x^{\frac{1}{2}} b(\x^{\frac{1}{2}})\\
\pm \log{\nu(\x)} &=& \kappa\tau\x^{-\frac{1}{2}}+\im 
\tilde{k}\pi+\x^{\frac{1}{2}} \tilde{b}(\x^{\frac{1}{2}})\\
\end{eqnarray*} where $uv= -\kappa^{2}dz^{2}$ is the quadratic invariant 
defined on page~\pageref{defn:uv}.
\end{prop}

With these propositions in mind, we define a smooth hyperelliptic 
curve $\hat{X}$ corresponding to each non-trivial solution $(A,\Phi)$ 
of \ref{eq:}.
Let 
$\alpha_{1},\ldots,\alpha_{m},\bar\alpha_{1}^{-1},\ldots,\bar\alpha_{m}^{-1}$ 
be the odd order zeros of $h(\x)^{2}-4$. 
\begin{enumerate}
\item If $\det\Phi\neq 0$, let $\hat{X}$ be 
the curve 
$$y^{2}=\prod_{i=1}^{m}(\x-\alpha_{i})(\x-\bar\alpha^{-1}).$$
\item If $\det\Phi = 0$, let $\hat{X}$ be 
the curve
$$y^{2}=\x\prod_{i=1}^{m}(\x-\alpha_{i})(\x-\bar\alpha^{-1}).$$
\end{enumerate}

Denote by $\pi:\hat{X}\ra\mathbb{C}P^{1}$ the projection map 
$(\x,y)\mapsto 
\x$. 

By definition, $\mu$ and $\nu$ are well-defined 
functions on $\hat{X}$ and 
$$\Theta:=\frac{1}{2\pi\im}d\log\mu\mbox{, 
}\Psi:=\frac{1}{2\pi\im}d\log\nu$$ are 
meromorphic differentials whose only singularities are double poles 
at $\pi^{-1}\{0,\infty\}$, and which have no residues. 

$\hat{X}-\pi^{-1}\{0,\infty\}$ is the Riemann surface of the 
eigenvalues of the holonomy matrices 
$H(\x,p)$. It supports a natural line bundle for each $p\in 
T^{2}$, namely the bundle 
of eigenspaces of $H(\x,p)$:
$$\left(E_{p}\right)_{(\x,y)}\subseteq 
\ker(H(\x,p)-\mu(\x,y)).$$
Since the fundamental group  $\pi_{1}(T^{2},p)$ 
is abelian, these eigenspaces are independent of the choice of 
generator. This line bundle can be naturally 
extended across $\pi^{-1}\{0,\infty\}$ to give a holomorphic line 
bundle $E_{p}$ on $\hat{X}$. 
\begin{prop}
Let $\hat{X}$ be the smooth hyperelliptic curve associated to a 
solution $(A,\Phi)$ of (\ref{eq:}). Then:
\begin{enumerate}
\item $\hat{X}$ has a real structure (i.e. an antiholomorphic 
involution)  $\rho:\hat{X}\ra\hat{X}$ that commutes with $\pi$ and 
covers the real structure 
$\x\mapsto\bar{\x}^{-1}$ of $\mathbb{C}P^{1}$. 
\item The hyperelliptic involution $\sigma:\hat{X}\ra\hat{X}$ 
commutes 
with $\pi$ and has no real fixed points.  
\item The differentials $\Theta$ and $\Psi$ satisfy
$$\sigma^{*}\Theta = -\Theta,\;\sigma^{*}\Psi = 
-\Psi,\;\rho^{*}\Theta = \bar\Theta,\;\rho^{*}\Psi = \bar\Psi.$$
\item $\Theta$ and $\Psi$ are linearly independent over $\mathbb{R}$.
\item For each $p\in T^2$, $j$ gives an isomorphism between $(\sigma\rho)^{*}E_p$ and $E_p$ 
whose square is $-1$, so each $E_p$ is quaternionic with respect to the 
real 
structure $\sigma\rho$. 
\item The periods of $\Theta$ and $\Psi$ are all integers.
\end{enumerate}
\end{prop}
The involution $\rho$ is the lift of $x\ra\bar{x}^{-1}$ that fixes the points in 
$\pi^{-1}\{x:|x|=1\}$. Using the fact that $d_{\x}$ is unitary 
for $\x$ on the unit circle,
$$\rho^{*}\mu=\bar{\mu}^{-1}$$ and by definition of $\hat{X}$,
$$\sigma^{*}\mu=\mu^{-1}.$$
The most important condition in the Proposition above is the last 
one, 
as the existence of meromorphic differentials with integral periods 
places a stringent restriction on $\hat{X}$.

The algebraic curve that we shall associate to a harmonic map $\map$ 
is a possibly singular curve $X$ (called the {\it spectral curve}) of 
which $\hat{X}$ is the 
normalisation. $X$ will reflect the geometry of the eigenspaces of 
the holonomy rather than merely that of the eigenvalues. We would 
like to be able to employ our algebro-geometric description to study 
families of harmonic maps, and since the limit of a family of 
smooth curves may be a singular one, it is sensible to allow spectral 
curves to be singular. The spectral curve will also enable us to 
compute the degree of the eigenspace bundle.  
Since $V$ has rank two, for $v,w\in V_{p}$ we have
$$\omega_{p}(v,w)=0 \iff\mbox{ v and w are linearly dependent.}$$
Thus $\omega_{p}$ vanishes precisely at those points $(\x,y)\in\hat{X}$ 
 at which
$\left(E_{p}\right)_{(\x,y)}$ and 
$\left(\sigma^{*}E_{p}\right)_{(\x,y)}$ 
coincide as subspaces of $V_{p}$. It necessarily vanishes at the 
branch points of $\hat{X}$, though for those branch points other 
than $0$ and $\infty$, it may do so to some odd order 
\mbox{$> 1$. }Suppose that $\omega$ vanishes to 
order $2k_{i}+1$ at $\alpha_{i}$ and $\bar\alpha_{i}^{-1}$, and to 
order $l_{j}$ at non-branch points $\beta_{j}$, $j=1,\ldots r$ of 
$\hat{X}$.

Then $X$ is described by the equation
$$y^{2}=\prod_{i=1}^{m}(\x-\alpha_{i})^{2k_{i}+1}(\x-\bar\alpha_{i}^{-1})^{2k_{i}+1}\prod_{j=1}^{r}(\x-\beta_{j})^{2l_{j}}.$$
Our definition of $X$ uses a particular eigenspace bundle $E_{p}$, 
but since for $\x\in\mathbb{C}^{*}$, the eigenspaces of the 
holonomy 
matrices $H(\x,p)$ and $H(\x,q)$ are related by parallel 
translation of $d_{\x}$, $X$ is in fact independent of the 
choice 
of base point $p$. It follows from the definition of $X$ that for 
each 
$p\in T^{2}$, $E_{p}$ is defined as a bundle over $X$. Furthermore, 
from 
the adjunction formula we may deduce that 
$$deg(E^{*}_{p})=g_{a}+1 \mbox{, where $g_{a}$ is the arithmetic 
genus of $X$.}$$
We may then define a map
$$\begin{array}{cccc}
l:&T^{2}&\ra&Pic^{g+1}X\\
&p&\mapsto&E_{p}^{*}
\end{array}$$

Given the real structure $\sigma\rho:X\ra X$, we say that a 
holomorphic line bundle $L$ on $X$ is {\it real} if there is an isomorphism 
$i:L\ra\overline{\sigma\rho^{*}L}$ whose square gives multiplication by a 
positive scalar on $L$. $L$ is {\it quaternionic} if there is an 
isomorphism 
$i:L\ra\overline{\sigma\rho^{*}L}$ whose square gives multiplication by a 
negative scalar on $L$. The line bundles 
$E_{p}^{*}$ are all quaternionic, and so if we fix $p\in T^{2}$, then 
for each $q\in T^{2}$, $ E_{q}\otimes E_{p}^{*}$ is real. In fact the 
map 
$$\begin{array}{cccc}
l\otimes E_{p}^{*}:&T^{2}&\ra&Pic^{0}X\\
&q&\mapsto&E_{q}\otimes E_{p}^{*}
\end{array}$$
is {\it linear}, and so is a linear map of $(T^{2},\tau)$ to a real 
torus in $Pic^{0}(X)$. If $g_{a}(X)\geq 2$, then this map is 
injective. There is a natural isomorphism 
$$(V)_{p}^{*}\cong H^{0}(X,E_{p}^{*}),$$
which enables us to reconstruct  $V$ from the eigenspace bundles.
The space of quaternionic line bundles of degree 
$g_{a}(X)$ is connected, and  the line bundles $E_{p}^{*}$ are 
non-special. 
The following is essentially Theorem 8.1 of \cite{hitchin:90}.
\begin{thm}\label{thm:general}
Let $X$ be a hyperelliptic curve $y^{2}=P(\x)$ of arithmetic 
genus $g_{a}$, and let \mbox{$\pi:X\ra\mathbb{C}P^{1}$} be the 
projection \mbox{$\pi(\x,y)=\x$.} Suppose $X$  satisfies:
\begin{enumerate}
\item $P(\x)$ is real with respect to the real structure 
\mbox{$\x\mapsto\bar{\x}^{-1}$} on $\mathbb{C}P^{1}$.
\item $P(\x)$ has no real zeros (i.e. no zeros on the unit 
circle $\x = \bar\x^{-1}$).
\item $P(\x)$ has at most simple zeros at $\x=0$ and 
$\x=\infty$.
\item There exist differentials $\Theta$ and $\Psi$ of the second 
kind 
on $X$ with periods in $\mathbb{Z}$.
\item $\Theta$ and $\Psi$ have double poles at $\pi^{-1}(0)$ and 
$\pi^{-1}(\infty)$ and satisfy
\begin{center}$\sigma^{*}\Theta=-\Theta,\,\sigma^{*}\Psi=-\Psi,\,
\rho^{*}\Theta=\bar\Theta,\,\rho^{*}\Psi=\bar\Psi$\end{center} 
where $\sigma$ is the hyperelliptic involution 
$(\x,y)\mapsto(\x,-y)$ and $\rho$ is the real structure 
induced from $\x\mapsto\bar\x^{-1}$.
\item The principal parts of $\Theta$ and $\Psi$ are linearly 
independent over $\mathbb{R}$.
\end{enumerate}
Then, for each point $E_0$ in the Picard variety of line bundles of degree 
$g_{a}+1$ on $X$ which are quaternionic with respect to the real 
structure $\sigma\rho$, there is a solution of (\ref{eq:}) for a 
torus, such that $X$ is the spectral curve of the solution and 
$\Theta$, $\Psi$ the corresponding differentials. The solution is, 
moreover, unique modulo gauge transformations and the operation of 
tensoring $V$ by a flat $\mathbb{Z}_{2}$-bundle. (Note that 
\setcounter{enumi}{2}({\it\arabic{enumi}}) implies that $\rho\sigma$ has 
no 
fixed points and so by \cite{Atiyah:71} quaternionic bundles of degree $g_{a}+1$ 
exist.) 
\end{thm}
Begining with the spectral data $(X,\Theta,\Psi,E_{0})$, one may construct a 
solution $(A,\Phi)$ to (\ref{eq:}) by carefully reversing the steps 
outlined above. For details, the reader is referred to \cite{hitchin:90}.

Given such data, we may calculate the 
eigenvalues of the holonomy using the relations
\begin{equation}\Theta=\frac{1}{2\pi\im}d\log\mu,\,\Psi=\frac{1}{2\pi\im}d\log\nu \label{eq:dlog}\end{equation}
on $X-\pi^{-1}\{0,\infty\}$. 
The solutions are unique only up to multiplication by a constant, but 
if we demand that \begin{equation} 
\mu\sigma^*\mu=1,\,\nu\sigma^*\nu=1,\label{eq:sigmamu}\end{equation} 
they are defined up to sign. 
By choosing $\mu,\,\nu$ satisfying  (\ref{eq:dlog}) and 
(\ref{eq:sigmamu}), we determine a solution $(A,\Phi)$ of (\ref{eq:}) 
up to gauge equivalence, that is, up to a diffeomorphism of the 
principal bundle $P$ that covers the identity on $T^2$ and commutes 
with the action of $SU(2)$. As mentioned earlier, this solution 
corresponds to a harmonic map $\map$ if and only if the flat 
connections $d_1$ and $d_{-1}$ are {\it trivial}. Since they are in any 
case unitary, this occurs precisely when $\mu$ and $\nu$ take the 
value 1 at all points in $\pi^{-1}\{ 1,-1\}$. Given such a solution $(A,\Phi)$, let $s_1$, $s_{-1}$ be  constant covariant sections of $P$ with respect to the connections $d_1$ and $d_{-1}$ respectively. We define $f:(T^2,\tau)\ra SU(2)$ by $s_1(p)=s_{-1}(p)f(p)$. Thus $f$ is unaffected by gauge transformations, but if we choose different covariant constant sections $\tilde{s}_1=s_1h$, $\tilde{s}_{-1}=s_{-1}k$, then 
$$\tilde{s}_1=\tilde{s}_{-1}k^{-1}fh,$$
so the map $f$ is well-defined modulo right and left actions of $SU(2)$, or modulo the action of $SO(4)=SU(2)\times SU(2)/\{\pm 1\}$ on $S^3$.

The following is Theorem 8.20 of \cite{hitchin:90}.
\begin{thm}\label{thm:ftns}
Let $(X,\Theta,\Psi,E_0)$ be spectral data satisfying the 
conditions of Theorem~\ref{thm:general}, where $X$ is given by $y^2=P(\x)$. Let $\mu$ and $\nu$ be 
functions on $X-\pi^{-1}\{0,\infty\}$ satisfying \newline 
$\Theta=\frac{1}{2\pi\im}d\log\mu,\,\Psi=\frac{1}{2\pi\im}d\log\nu $ and 
$\mu\sigma^*\mu=1,\,\nu\sigma^*\nu=1$. Then
\begin{enumerate}
\item $(X,\mu,\nu)$ determines a harmonic map from a torus to $S^3$ 
if and only if $$\mu(\x,y)=\nu(\x,y)=1\mbox{ for all 
}(\x,y)\in\pi^{-1}\{1,-1\}.$$
\item The map is conformal if and only if $P(0)=0$.
\item The torus maps to a totally geodesic 2-sphere if and only if 
$g_a$ is odd, $P(\x)$ is an even polynomial, and the point 
$E_0\in{\rm Pic}^{g_a+1}(X)$ and the functions $\mu$ and $\nu$ on 
$X-\pi^{-1}\{0,\infty\}$ are invariant under $\sigma\gamma$, where 
$\gamma$ is the involution of $X$ defined by 
$\gamma(\x,y)=(-\x,y)$.
\item The harmonic map is uniquely determined by $(X,\mu,\nu, E_0)$ modulo 
the action of $SO(4)$ on $S^3$.
\end{enumerate}
\end{thm}

\chapter{Conformal Maps}
We prove the following theorem:
\begin{thm}
For each integer $g>0$ there are countably many conformal harmonic immersions from rectangular tori to $S^3$
whose spectral curves have genus $g$.
\end{thm}

We in fact demonstrate  the existence of spectral curves possessing an
additional symmetry, namely \mbox{$x\mapsto \frac{1}{x}$}, where
\proj. This symmetry induces two holomorphic involutions on $X$,
which we utilise by quotienting out by them. The resulting quotient curves
$\Cpm$ are our basic object of study, and we obtain differentials on
$X$ by pulling back differentials from \mbox{$C_{\pm}$}. We will use
proof by induction, in which at each
induction step the genera of $\Cpm$ increase by one, and hence the genus $g$ of $X$ increases by two. Thus we divide our proof into the even and odd genus cases.
Our proofs extend the methods of Ercolani, Kn\"{o}rrer and Trubowitz
~\cite{EKT:93}, who showed that for each even genus $g\geq 2$ there
is a constant mean curvature torus in $\mathbb{R}^3$ whose spectral curve has
genus $g$. (There is also a spectral curve construction for constant mean curvature tori.)

\section{Odd Genera}
Let $C_+=C_+(R,\lambda_1,\bar\lambda_1,\ldots,\lambda_n,\bar\lambda_n)$ be the curve given by
$$w_{+}^2=(z-R)\prod_{i=1}^{n}(z-\lambda_i)(z-\bar{\lambda_i})$$ and
$C_-=C_-(R,\lambda_1,\bar\lambda_1,\ldots,\lambda_n,\bar\lambda_n)$ that given by
$$w_{-}^2=(z-2)(z+2)(z-R)\prod_{i=1}^{n}(z-\lambda_i)(z-\bar{\lambda_i})$$
where we assume that \mbox{$R\in (-\infty,-2)\cup(2,\infty)$}, \mbox{$\lambda_i\neq 2$} for \mbox{$i=1,\ldots,2n$} and
\mbox{$\lambda_i\neq\lambda_j$} for $i\neq j$. Let
\mbox{$\pi_{\pm}:\Cpm\ra\mathbb{C}P^1$}, $(z,w_\pm)\mapsto z$ denote
the respective projections to the Riemann sphere. (We shall
henceforth generally
omit the word ``respective''.) Construct
\mbox{$\pi:X\ra\mathbb{CP}^1$} as the fibre product of
\mbox{$\pi_{+}:C_+\ra\mathbb{C}P^1$} and
\mbox{$\pi_{-}:C_-\ra\mathbb{C}P^1$}, that is, let
\mbox{$X:=\{(p_+,p_-)\in C_+\times C_-: \pi_+(p_+)=\pi_-(p_-)\}$},
with the induced algebraic structure and the obvious projection
$\pi$ to \cp. Notice then that $X$ is given by the equation
$$y^2=x(x-r)(x-r^{-1})\prod_{i=1}^{n}(x-\alpha_i)(x-{\alpha_i}^{-1})(x-\bar{\alpha_i})(x-{\bar{\alpha_i}}^{-1}),$$
where $$r+r^{-1}=R,\,\alpha_i+{\alpha_i}^{-1}=\lambda_i,$$
and $\pi$ is given by
$$\pi:(x,y)\mapsto x.$$ These identifications occur via the maps
$$q_{+}(x,y)=\left(x+\frac{1}{x},\frac{y}{x^{n+1}}\right)=(z,w_+)$$
and
$$q_{-}(x,y)=\left(x+\frac{1}{x},\frac{(x+1)(x-1)y}{x^{n+2}}\right)=(z,w_-).$$

$X$ has genus $2n+1$ and possesses the holomorphic involutions
$$\begin{array}{rccc}
i_{\pm}:&X&\longrightarrow&X\\
&(x,y)&\longmapsto&\left(\frac{1}{x},\frac{\pm y}{x^{2n+1}}\right).
\end{array}$$
The curves $\Cpm$ are the quotients of $X$ by these involutions, with
quotient maps $q_\pm:X\ra\Cpm$.
$\Cpm$ each possesses a real structure $\rho_\pm$, characterised by the
properties that it covers the involution $z\mapsto\bar z$ of
$\mathbb{C}P^1$ and fixes the points in
$\pi_\pm\!^{-1}[-2,2]$. These real structures are given by
$$\rho_\pm(z,w_\pm)=(\bar z,\mp\bar w_\pm),\mbox{ for $R>2$}$$ or
$$\rho_\pm(z,w_\pm)=(\bar z,\pm\bar w_\pm),\mbox{ for $R<-2$.}$$
The cases $R>2$ and $R<-2$ are similar, but the sign difference
carries through to future computations. For simplicity of exposition
we assume henceforth that $R>2$. Then the corresponding real
structure on $X$ is given by
$$\rho(x,y)=\left(\frac{1}{\bar x},\frac{-\bar y}{{\bar x}^{2n+1}}\right).$$

Whilst our primary interest lies with curves $C_\pm$ as described above, we
 consider also curves \mbox{$\Cpm=\Cpm(R,\lambda_1,\ldots,\lambda_{2n})$} given by

$$w_{+}^2=(z-R)\prod_{i=1}^{2n}(z-\lambda_i)$$ and

$$w_{-}^2=(z-2)(z+2)(z-R)\prod_{i=1}^{2n}(z-\lambda_i)$$ respectively,
where we assume that \mbox{$R\in (2,\infty)$}, \mbox{$\lambda_i\neq\pm 2$} for \mbox{$i=1,\ldots,2n$} and that the sets
\mbox{$\{\lambda_1,\lambda_{2}\},\ldots,\{\lambda_{2n-1},\lambda_{2n}\}$} are mutually disjoint.

Take $(R,\lambda_1,\ldots,\lambda_{2n})$ as described above. Let \mbox{$\tilde a_0,\ldots,\tilde a_n $} be simple closed curves in
\mbox{$\cp-\{R,2,-2,\lambda_1,\ldots,\lambda_{2n}\}$}, and $\tilde c_1, \tilde c_{-1}$ simple closed curves in \mbox{$\cp-\{R,\lambda_1,\ldots,\lambda_{2n}\}$}, such that
\begin{enumerate}
\item $\tilde a_0$ has winding number one around $2$ and $R$, and winding number zero around the other branch points of $C_-$,
\item for $i=1\ldots n$, $\tilde a_i$ has winding number one around $\lambda_{2i-1}$ and $ \lambda_{2i}$, and winding number zero around the other branch points of $C_\pm$,
\item $\tilde c_1$ begins and ends at $z=2$, has winding number one around $R$ and zero around $\lambda_i$, $i=1,\ldots,2n$,
\item $\tilde c_{-1}$ begins and ends at $z=-2$, and has winding number one around $R$ and each $\lambda_i$ $i=1,\ldots,2n$.
\end{enumerate}

Choose lifts of the curves \mbox{$\tilde a_1,\ldots,\tilde a_n$} to  $C_+$ and also of \mbox{$\tilde a_0,\ldots,\tilde a_n$} to $C_-$. Let \mbox{$a^-_0,a^\pm_1,\ldots,a^\pm_n\in H_1(C_{\pm},\mathbb{Z})$} denote the homology classes of these lifts. Denote by \mbox{$b^-_0,b^\pm_1,\ldots,b^\pm_n$} the
completions to canonical bases of  \mbox{$H_1(\Cpm,\mathbb{Z})$.} Choose open curves $c_1$, $c_{-1}$  in $C_+$ covering the loops
$\tilde c_1$ and $\tilde c_{-1}$.

\begin{figure}[h]
\hspace*{1.5cm}\includegraphics{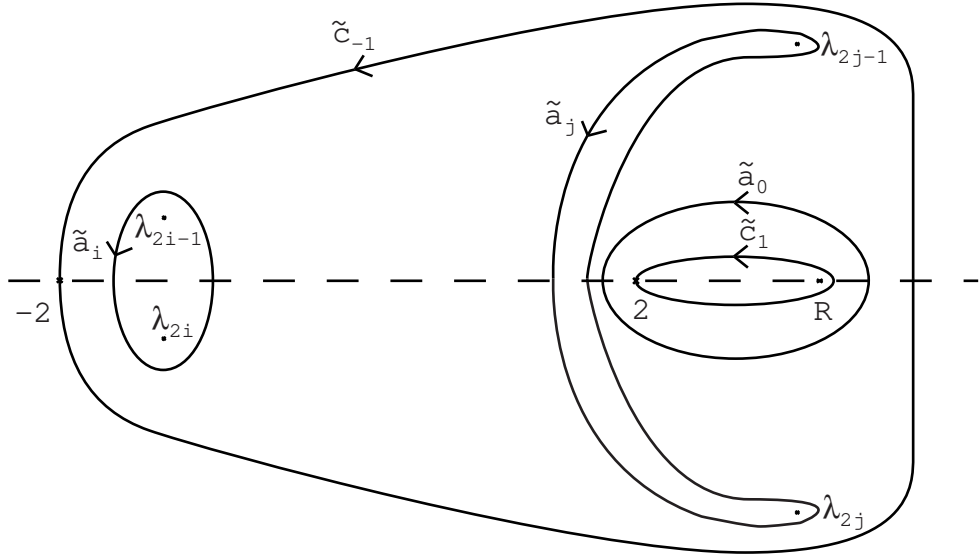}
\caption{Curves $\tilde a_i$, $\tilde c_{\pm 1}$ for $p\in M_{n,\R}$.\label{fig:Cpm}}
\end{figure}

Denote by $M_n$ the space of $2n+1$-tuples $(R,\lambda_1,\ldots,\lambda_{2n})$ as above
together with the choices we have described. Let $M_{n,\R}$ denote the subset of $M_n$ such that (see figure~\ref{fig:Cpm}):
\begin{enumerate}
\item $\lambda_{2i}=\bar\lambda_{2i-1}$ for $i=1,\ldots,n$,
\item for $i=1\ldots n$, $\tilde a_i$ is invariant under conjugation and intersects the real axis exactly twice, both times in the interval \mbox{$(-2,2)$},
\item the lifts of  \mbox{$\tilde
a_1,\ldots,\tilde a_n$} to $C_{+}$  are chosen so
that the point where $\tilde a_i$ intersects the $z$-axis with positive orientation is lifted to a point where $\frac{w_+}{\im}$ is negative,
\item the lifts of \mbox{$\tilde a_0,\ldots,\tilde a_n$} to $C_{-}$  are chosen so
that the point
where $\tilde a_i$ intersects the $z$-axis with positive orientation
is lifted to a point in $C_-$ where $w_-$ is positive,
\item $c_1$, $c_{-1}$ begin at points with
$\frac{w_+}{\im}<0$.
\end{enumerate}

\begin{figure}[h]\hspace*{1.2cm}\includegraphics{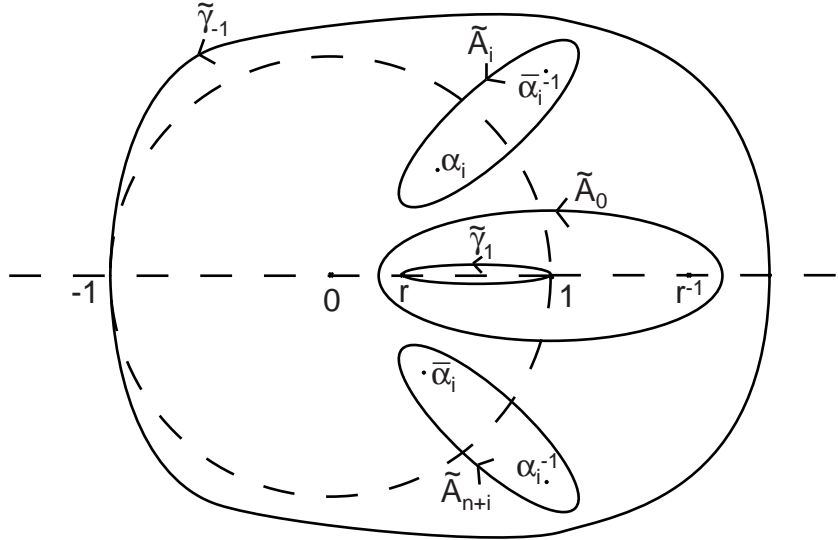}
\caption{The curves $\tilde A_i$ and $\tilde\gamma_{\pm 1}$.\label{fig:X}}
\end{figure}

For each $p\in M_{n,\R}$ there is a unique canonical basis
\mbox{$A_0,\ldots,A_{2n}$}, \mbox{$B_0,\ldots,B_{2n}$} for the
homology of $X$ such that \mbox{$A_0,\ldots,A_{2n}$} cover the
homotopy classes of loops \mbox{$\tilde A_0,\ldots,\tilde A_{2n}$}
shown in Figure~\ref{fig:X} and
$${(q_-)}_*(A_0)=2a^-_0,\;(q_\pm)_*(A_i)=\mp(q_\pm)_*(A_{n+i})=a^\pm_i.$$

There are also unique curves $\gamma_1$ and $\gamma_{-1}$ on $X$ such
that $(q_+)_*(\gamma_{\pm 1})=c_{\pm 1}$; they project to
$\tilde{\gamma_1}$ and $\tilde{\gamma_{-1}}$ of Figure~\ref{fig:X}. Note that
$\gamma_1$ connects the two points of $X$ with $x=1$ whilst
$\gamma_{-1}$ connects the two points of $X$ with $x=-1$.

Denote by ${\cal A}^\pm$ the subgroups of \mbox{$H_1(\Cpm,\mathbb{Z})$} generated by the $a^\pm$
classes. Then modulo ${\cal A}^\pm$,\nopagebreak
$${(q_-)}_*(B_0)\equiv b^-_0,\:(q_\pm)_*(B_i)\equiv
b^\pm_i,\:(q_\pm)_*(B_{n+i})\equiv \mp
b^\pm_i.$$

For each $p\in M_n$, we will define differentials $\Wpm(p)$ on $\Cpm(p)$, and show that for
each $n$ there exists \mbox{$p\in M_{n,\R}$}  such that appropriate integer multiples of ${q_+\!}^*(\Omega_+(p))$ and ${q_-\!}^*(\Omega_-(p))$ satisfy the conditions of Hitchin's correspondence. It is easier to argue this way than to perform a similar argument directly on $X$, for reasons that will become clear.

Let $p\in M_n$ and define $\Wpm=\Wpm(p)$ on $\Cpm(p)$ by:
\begin{enumerate}
\item $\Wpm(p)$ are meromorphic differentials of the second kind: their
only singularities are double poles at $z=\infty$, and they have no
residues.
\item $\int_{a^-_0}\Omega_-(p) = 0$ and $\int_{a^\pm_i}\Wpm(p) = 0$ for
$i=1,\ldots,n$.
\item As $z\ra\infty$, $\Omega_+(p)\ra\frac{z^n}{w_+(p)}$ and
$\Omega_-(p)\ra\frac{z^{n+1}}{w_-(p)}$.
\end{enumerate}
(We shall generally denote the paths of integration $a^\pm_i(p),b^\pm_i(p),c_{\pm}(p)$ simply by $a^\pm_i,b^\pm_i$ and $c_{\pm}$, as the point $p$ in question is usually clear from the differential being integrated.) In view of the defining conditions above, we may write
$$\Wp=\frac{\prod_{j=1}^n(z-\zeta^+_j)dz}{w_+}$$
and
$$\Wm=\frac{\prod_{j=0}^n(z-\zeta^-_j)dz}{w_-}.$$

Define
$$I_{+}(p):=\im\left(\int_{c_1}\Omega_+(p),\int_{c_{-1}}\Omega_+(p),\int_{b^+_1}\Omega_+(p),\ldots,\int_{b^+_n}\Omega_+(p)\right),$$

$$I_{-}(p):=\left(\int_{b^-_0}\Omega_-(p),\int_{b^-_1}\Omega_-(p),\ldots,\int_{b^-_n}\Omega_-(p)\right).$$

Then  for $p\in M_{n,\R}$, $I_+(p)$ and $I_-(p)$ are real, since then $$(\rho_\pm)_*(b^\pm_i)=b^\pm_i
\mbox{ mod ${\cal A}^\pm$},\:  (\rho_+)_*(c_{\pm 1})=c_{\pm 1} \mbox{
mod ${\cal A}^+$},$$ so
$$\int_{b^+_i}\Omega_+(p)=\int_{(\rho_+)_*(b^+_i)}\Omega_+(p)=\int_{b^+_i}\rho^*_+(\Omega_+(p))=-\int_{b^+_i}\overline\Omega_+(p),\:i=1\ldots
n,$$
and similarly $$\int_{c_{\pm 1}}\Omega_+(p)=-\int_{c_{\pm
1}}\overline\Omega_+(p)$$ and
$$\int_{b^-_i}\Omega_-(p)=\int_{b^-_i}\overline\Omega_-(p),\:i=0\ldots n.$$

Given $p\in M_{n,\R}$, there are real numbers $s_+$ and $s_-$ such that
$\im s_+q_+\!^*(\Omega_+(p))$ and $s_-q_-\!^*(\Omega_-(p))$ are differentials
on $X$ satisfying the conditions of Hitchin's correspondence if and
only if $I_+(p)$ and $I_-(p)$ represent rational elements of
$\mathbb{R}P^{n+1}$ and $\mathbb{R}P^{n}$ respectively.

\begin{thm}For each non-negative integer $n$, there exists \mbox{$p\in M_{n,\R}$} such that
\begin{enumerate}
\item $\zeta^+_j(p),\,j=1,\ldots,n$ are pairwise distinct, as are $\zeta^-_j(p),\, j=0,\ldots,n$.
\item The map \vspace*{-5mm}$$\begin{array}{rccc}\\
\phi:&M_n&\longrightarrow&\mathbb{C}P^{n+1}\times\mathbb{C}P^n\\
&p&\longmapsto&([I_+(p)],\,[I_-(p)]).
\end{array}$$
has invertible differential at $p$.
\end{enumerate}
\label{thm:odd}
\end{thm}
This gives that the restriction
$$\phi\left.\right|_{M_{n,\R}}:M_{n,\R}\ra\mathbb{R}P^{n+1}\times\mathbb{R}P^n$$
of $\phi$ to $M_{n,\R}$ also has invertible differential at $p$. Since rationality is a dense condition, the Inverse Function Theorem
then implies that for each positive odd integer $g$, there are
countably many spectral curves $X$ of genus $g$ each giving rise to a
torus $(T^2,\tau)$ and a branched minimal immersion \map. The
conformal type of the torus is given by

$$\tau = \frac{\im s_+\rm{p.p.}_\infty(q_+\!^*(\Omega_+(p)))}{s_-\rm{p.p.}_\infty(q_-\!^*(\Omega_-(p)))}=\frac{\im s_+}{s_-},$$
where $\rm{p.p.}_\infty(q_\pm\!^*(\Omega_\pm))$ denotes the principal
part of $q_\pm\!^*(\Wpm)$ at $\infty$. Thus each torus $(T^2,\tau)$
is rectangular.

In fact we shall prove a slightly stronger result, the extra strength
residing in a statement that arises naturally from an attempt to prove
Theorem~\ref{thm:odd} by induction on $n$, and enables one to
complete the induction step. This statement is somewhat lengthy to formulate, and will appear unmotivated at this juncture. Our approach is thus to present an attempt to prove Theorem~\ref{thm:odd} by induction, and derive the necessary modifications. The reader who wishes to view the modified statement at this point is referred to page~\pageref{thm:oddextra}.

\noindent{\bf Proof of Theorem~\ref{thm:odd}:}
We use induction upon $n$. We shall begin with the induction step, in order to formulate the
``extra conditions'' mentioned above. Suppose then that $\phi$ has invertible differential at \mbox{$p\in M_{n,\R}$}. For \mbox{$\m\in (-2,2)$}, \mbox{$\n\in\mathbb{R}$}, we shall denote by $(p,\m,\n)$ the point of \mbox{$M_{n+1,\R}$} such that $$\lambda_i(p,\m,\n)=\lambda_i(p),\,i=1,\ldots n,\; \lambda_{n+1}(p,\m,\n)=\m+\im\n.$$
Denote by \mbox{$(p,\m,\n)$} the point in \mbox{$M_{n+1,\R}$} with branch points $$\lambda_i(p,\mu,\nu)=\lambda_i(p),\,i=1,\ldots,2n$$ and
$$\lambda_{2n+1}(p,\mu,\nu)=\mu+\im\nu,\;\lambda_{2n+2}(p,\mu,\nu)=\mu-\im\nu.$$
We wish to show that for a  generic \mbox{$\m\in (-2,2)$} and $\n$ sufficiently small, $\phi$ has invertible differential at \mbox{$(p,\m,\n)$}. (Here ``generic'' means ``outside the zero set of a real-analytic
function''.) Once we have incorporated our ``additional condition'', we shall achieve our aim by considering the boundary case $\n=0$. Notice that choosing \mbox{$\pi_\pm :C_{\pm}\ra\mathbb{C}P^{1}$} to each  have an additional
branch point of multiplicity two on the interval $(-2,2)$ corresponds to
choosing \mbox{$\pi:X\ra\mathbb{C}P^{1}$} to have an additional branch point of multiplicity two on the unit circle. For brevity of notation, we shall write \mbox{$\po=\po(p,\m)$} for $(p,\m,0)$.

\begin{figure}[h]
\hspace*{1.5cm}\includegraphics{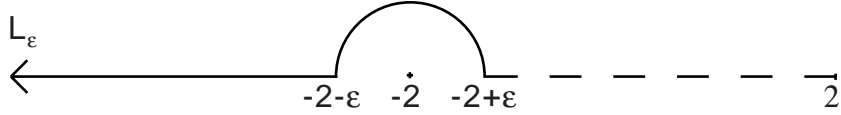}
\caption{$h_\epsilon$ is a function of $\m\in L_\epsilon$.\label{fig:hdomain}}\end{figure}

Let
\be H(\m):=\left(\begin{array}{cc}
I_{+}(\po)&0\\
0&I_{-}(\po)\\
\frac{\partial}{\partial R} I_{+}(\po)&\frac{\partial}{\partial R} I_{-}(\po)\\
\frac{\partial}{\partial \lambda_{1}}I_{+}(\po)&\frac{\partial}{\partial \lambda_{1}} I_{-}(\po)\\
\vdots&\vdots\\
\frac{\partial}{\partial \lambda_{2n}}I_{+}(\po)&\frac{\partial}{\partial \lambda_{2n}}I_{-}(\po)\\
\frac{\partial}{\partial \m}
I_{+}(\po)&\frac{\partial}{\partial \m} I_{-}(\po)\\
\frac{\partial^{2}}{\partial {\n}^{2}}
I_{+}(\po)&\frac{\partial^{2}}{\partial {\n}^{2}} I_{-}(\po)
\end{array}\right),\label{eq:h}\ee
and $$h(\m):=\det H(\m).$$
$h$ is a  real-analytic function of $\m\in(-2,2)$ and for each
\mbox{$\epsilon\in(0,\min_{i=1,\ldots, n}|\lambda_i+2|)$} we may use the
above
formula to define it as a real-analytic function $h_\epsilon$ of $\m$ on the curve
$L_\epsilon$ shown in figure~\ref{fig:hdomain}.

We will show that $h(\m)\not\ra 0$ as $\m\ra\infty$ along  $L_\epsilon$,
by computing asymptotics for each of the vectors in $H(\m)$. The advantage of passing to the quotient curves $\Cpm$ is the availability of this limiting argument. The interval $(-2,2)$ possesses a natural extension to a line whereas the unit circle does not. We will also prove that
$$\frac{\partial}{\partial
\n€}\left(I_{+}(\po),I_{-}(\po)\right)=0.$$

Then for generic $\mu\in(-2,2)$, we will have  that $h(\m)\neq 0$, and utilising
\be \frac{\partial}{\partial
{\n}}\left(I_{+}(p,\m,\n),I_{-}(p,\m,\n)\right)=
 \frac{\partial^{2}}{\partial {\n}^{2}}
\left(I_{+}(\po),I_{-}(\po)\right) +
O(\n^{2}),\label{eq:nusmall}\ee
this gives that for $\n$ sufficiently small, $d\phi_{(p,\m,\n)}$ is invertible.

One simplification provided by choosing $\n=0$ is that $C_\pm(p)$ are the respective normalisations of $C_\pm(\po)$, with normalisation maps

\begin{center}
$\begin{array}{rccc}
\Psi_\pm:&C_\pm(p)&\longrightarrow&C_\pm(\po)\\
&(z,w_\pm(p))&\longmapsto&(z,(z-\m)w_\pm(p))
\end{array}$
\end{center}
and that
$$\Psi_\pm\!^*(\Wpm(\po))=\Wpm(p),$$
whilst
\begin{eqnarray*}
(\Psi_+)_*(b^+_i(p))&=&b^+_i(\po), i=1,\ldots,n,\\
(\Psi_+)_*(c_{\pm 1}(p))&=&c_{\pm 1}(\po),\\
(\Psi_-)_*(b^-_i(p))&=&b^-_i(\po), i=0\ldots n.
\end{eqnarray*}

Thus
\begin{eqnarray}
I_+(\po)&=&\left(I_+(p),\im\int_{b^+_{n+1}}\Omega_+(\po)\right), \label{eq:Iplus}\\
I_-(\po)&=&\left(I_-(p),\int_{b^-_{n+1}}\Omega_-(\po)\right).\label{eq:Iminus}
\end{eqnarray}
For each point $q\in M_n$, let $u_\pm(q)$ be local coordinates on $C_{\pm}(q)$ near
$\pi_{\pm}\!^{-1}(\infty)$ such that
$$u_{\pm}(q)^{2}=z_{\pm}\!^{-1}$$
and
$$w_{+}(q)=u_{+}(q)^{2n+1}+O(u_{+}(q)^{2n}) \mbox{ as } z_{+}\ra\infty$$ whilst $$w_{-}(q)=u_-(q)^{2n+3}+O(u_{-}(q)^{2n+2}) \mbox{ as }z_{-}\ra\infty.$$
Then for $z$ near $\infty$,
\be \Wpm(q)=(u_{\pm}(q)+D_{\pm}(q)u_{\pm}(q)^{3}+O(u_{\pm}(q)^{5}))dz,
\label{eq:Dplusdefn}\ee
where
\be D_{+}(q):=\frac{1}{2}R(q)+\sum_{i=1}^{2n}\lambda_{i}(q)-\sum_{j=1}^n\zeta^+_j(q),\label{eq:Dplus}\ee and
\be D_{-}(q):=\frac{1}{2}R(q)+\sum_{i=1}^{2n}\lambda_{i}(q)-\sum_{j=0}^n\zeta^-_j(q).\label{eq:Dminus}\ee
We assumed that the $\zeta^+_j$ are pairwise distinct, $j=1,\ldots,n$.
Thus the differentials $\frac{\Wp(p)}{z-\zeta^+_j}$ are a basis for the holomorphic differentials on $C_+(p)$, and we may define
\mbox{$c^+_j(p)$} by the equations
\be \frac{3}{2}\int_{a^+_j}z\Wp(p)+\sum_{j=1}^{n}c^+_j(p)\int_{a^+_i}\frac{\Wp(p)}{z-\zeta^+_j}=0,\:i,j=1,\ldots,n,\label{eq:cj}\ee
and let
\be\Wph(p):=\frac{3}{2}z\Wp(p)+\sum_{j=1}^{n}c^+_j(p)\frac{\Wp(p)}{z-\zeta^+_j}\label{eq:Wph},\ee
$$\widehat{I}_+(p):=\im\left(\int_{c_1}\widehat{\Omega}_+(p),\int_{c_{-1}}\widehat{\Omega}_+(p),\int_{b^+_1}\widehat{\Omega}_+(p),\ldots,\int_{b^+_n}\widehat{\Omega}_+(p)\right),$$
$$\widehat{I}_-(p):=\left(\int_{b^-_0}\widehat{\Omega}_-(p),\int_{b^-_1}\widehat{\Omega}_-(p),\ldots,\int_{b^-_n}\widehat{\Omega}_-(p)\right).$$
\begin{lemma}
As $\m\ra\infty$ along $L_\epsilon$, the following asymptotic expressions hold:
\begin{enumerate}
\item
${\displaystyle I_{+}(\po)=\left(I_{+}(p),4\im\m^{1/2}-4\im
D_{+}(p)\m^{-1/2}+O(\m^{-3/2})\right)}$\\
${\ds I_{-}(\po)=\left(I_{-}(p),4\m^{1/2}-4D_{-}(p)\m^{-1/2}+O(\m^{-3/2})\right)}$

\item
${\ds\frac{\partial}{\partial{R}} I_+(\po)=\left(\frac{\partial}{\partial R} I_{+}(p),\im\left(-2+\sum_{j=1}^{n}\frac{\partial\zeta^{+}_{j}}{\partial
R}\right)\m^{-1/2}+O(\m^{-3/2})\right)}$\\
${\ds\frac{\partial}{\partial{R}} I_{-}(\po)=\left(\frac{\partial}{\partial R}
I_{-}(p),\left(-2+\sum_{j=0}^{n}\frac{\partial\zeta^{-}_{j}}{\partial
R}\right)\m^{-1/2}+O(\m^{-3/2})\right)}$

\item
For $ i=1,\ldots,2n$,\\
${\ds\frac{\partial}{\partial{\lambda_{i}}} I_+(\po)=\left(\frac{\partial }{\partial \lambda_{i}}
I_{+}(p),\im\left(-2+\sum_{j=1}^{n}\frac{\partial\zeta^{+}_{j}}{\partial
\lambda_{i}}\right)\m^{-1/2}+O(\m^{-3/2})\right)}$, and \\
${\ds\frac{\partial}{\partial{\lambda_{i}}} I_-(\po)=\left(\frac{\partial }{\partial \lambda_{i}}
I_{-}(p),\left(-2+\sum_{j=0}^{n}\frac{\partial\zeta^{-}_{j}}{\partial
\lambda_{i}}\right)\m^{-1/2}+O(\m^{-3/2})\right)}$

\item
${\ds\frac{\partial }{\partial\m}I_+(\po)=\left(0,2\im\m^{-1/2}+2\im
D_{+}(p)\m^{-3/2}+O(\m^{-5/2})\right)}$\vspace{3mm}\\
${\ds\frac{\partial}{\partial\m}
I_-(\po)=\left(0,2\m^{-1/2}+2D_{-}(p)\m^{-3/2}+O(\m^{-5/2})\right)}$

\item
${\ds\frac{\partial}{\partial\n}I_\pm(\po)=0}$

\item %
${\ds\frac{\partial^{2}}{\partial\n^{2}}I_+(\po)=\left(\left(-\frac{1}{2}\m^{-2}+D_+(p)\m^{-3}\right)I_{+}(p)-\m^{-3}\widehat{I}_{+}(p)+O(\m^{-4}),\right.}$
\vspace{-6mm}\begin{flushright}${\ds\left.\hfill \frac{3\im}{2}\m^{-3/2}+\frac{9\im}{2}D_{+}(p)\m^{-5/2}+O(\m^{-7/2})\right)}$\end{flushright}

${\ds\frac{\partial^{2}}{\partial\n^{2}}I_-(\po)=\left(\left(-\frac{1}{2}\m^{-2}+D_-(p)\m^{-3}\right)I_{-}(p)-\m^{-3}\widehat{I}_{-}(p)+O(\m^{-4}),\right.}$
\vspace{-6mm}\begin{flushright}${\ds\left.\frac{3}{2}\m^{-3/2}+\frac{9}{2}D_{-}(p)\m^{-5/2}+O(\m^{-7/2})\right)}$\end{flushright}
\end{enumerate}
\label{thm:asym}
\end{lemma}

\begin{figure}[h]\hspace*{1.5cm}\includegraphics{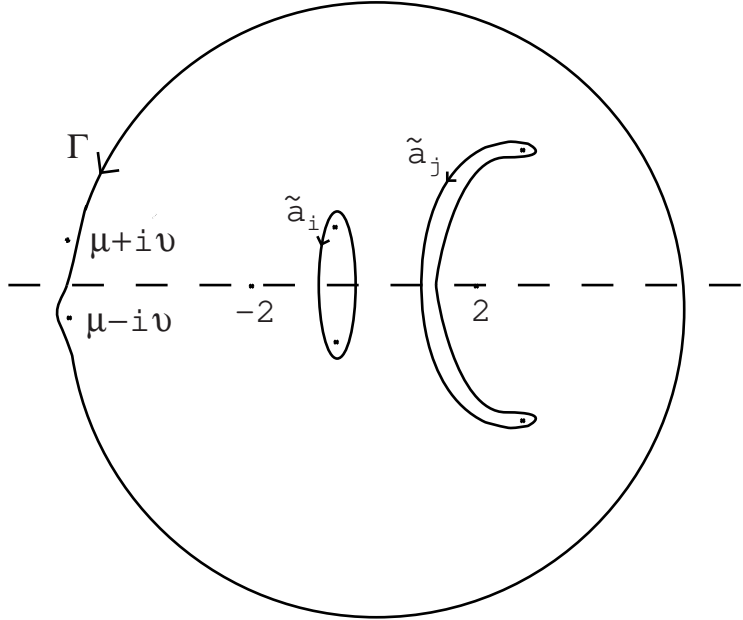}
\caption{We can take a representative of $b^+_{n+1}(\po)$ that projects to a circle.\label{fig:bnplusone}}\end{figure}

\noindent {\bf Proof of Lemma~\ref{thm:asym}:}
All but the last components of \setcounter{enumi}{1}({\it\arabic{enumi}})--\setcounter{enumi}{4}({\it\arabic{enumi}}) are applications
of equations~\ref{eq:Iplus} and \ref{eq:Iminus}. The last components
of \setcounter{enumi}{1}({\it\arabic{enumi}})--\setcounter{enumi}{4}({\it\arabic{enumi}}) involve integrals over the curves
$b^{\pm}\!_{n+1}(\po)$. Let $\Gamma$ denote the circle $|z|=\m$, transversed clockwise.  For $\m$ sufficiently large,
\begin{eqnarray*}\int_{b^+_{n+1}(\po)}\Wp(\po)
&=&-\int_{\Gamma}\frac{\prod_{j=1}^n(z-\zeta^+_j)dz}{\sqrt{(z-R)\prod_{i=1}^n (z-\lambda_i)(z-\bar\lambda_i)}}\mbox{ (see figure~\ref{fig:bnplusone})}\\
&=&4\m^{1/2}-4D_\pm(p)\m^{-1/2}+O(\m^{-3/2}),
\end{eqnarray*}
and similarly for \mbox{$\int_{b^-_{n+1}(\po)}\Omega_-(\po)$}, which easily gives the remainder of \setcounter{enumi}{1}({\it\arabic{enumi}})--\setcounter{enumi}{4}({\it\arabic{enumi}}).

\noindent{\bf\setcounter{enumi}{5}({\it\arabic{enumi}}):} We employ the notation $\dot f$ to indicate $\left.\frac{\partial f}{\partial\n}\right|_{\n=0}$. We will work only with $\Omega_+(\po)$, but similar arguments apply to $\Omega_-(\po)$. For $i=1,\ldots,n+1$,
\begin{equation}\int_{a^+_i}\dot\Omega_+(\po)=\left.\pnu\right|_{\n=0}\int_{a^+_i}\Omega_+(\po)=0.\label{eq:dotW}\end{equation}
For $\n$ small, we may write $$\Omega_+(p,\m,\n)=\frac{\prod_{j=1}^{n+1}(z-\zeta^+_j(p,\m,\n))dz}{w_+(p,\m,\n)},$$ where $\zeta^+_j$ are analytic functions satisfying $$\zeta^+_j(\po)=\left\{\begin{array}{ll}\zeta_j(p)&\mbox{ for }j=1,\ldots,n\\
\m&\mbox{ for }j=n+1.
\end{array}\right.$$ (Due to this accordance, we write simply $\zeta^\pm_j$ for $\zeta^\pm_j(p)$ or $\zeta^\pm_j(\po)$.)\\
Then
$$\dot\Omega_+(\po)=\sum_{j=1}^{n+1}\left(\frac{-\dot\zeta^+_j}{z-\zeta^+_j}\right)\frac{\prod_{k=1}^{n+1}(z-\zeta_k)dz}{w_+}.$$
If $\m$ coincides with one of the $\zeta_j$, $j=1,\ldots,n$, then we immediately see that the lift of $\dot\Omega_+$ to the normalisation $C_+(p)$ of $C_+(\po)$ is holomorphic. If $\m$ does not equal any of the $\zeta^+_j$, then we may again conclude that $\dot\Omega_+(\po)$ lifts to a holomorphic differential by employing the observation that
(\ref{eq:dotW}) for $j=n+1$ states that $\dot\Omega_+(\po)$ has zero residue at $z=\m$, and hence $\dot\zeta_{n+1}=0$. Equation~(\ref{eq:dotW}) tells us that this lift has zero $a$-periods, and so is itself zero.
For $i=1,\ldots,n$ then,
$$\left.\pnu\right|_{\n=0}\int_{b^+_i}\Omega_+(\po)=\int_{b^+_i}\dot\Omega_+(\po)=0.$$

To compute $\left.\pnu\right|_{\n=0}\int_{b^+_{n+1}}\Omega_+(\po)$ we use reciprocity with the holomorphic differential $\omega(p,\m,\n)$ on $C_+(p,\m,\n)$ defined by
$\int_{a^+_i}\omega(p,\m,\n) = \Bigg\{\begin{array}{l}0\mbox{, for $i=1,\ldots,n$}\\
2\pi\im\mbox{, for $i=n+1$.}\end{array}$\\
To simplify our notation, we write $\ppr$ for  $(p,\m,\n)$.
 Let $\Delta(\ppr)$ denote the polygon formed by cutting $C_{+}(\ppr)$ open
along representatives of the homology elements $a^+_{i}, b^+_{i}$, and
choose $q_{0}\in\Delta(\ppr)$. Define a holomorphic function $g_{\ppr}$ on
$\Delta(\ppr)$ by $$g_{\ppr}(q):=\int_{q_0}^{q}\Omega_+(\ppr).$$
Reciprocity gives
$$2\pi\im\sum_{\rm residues}g_{\ppr}\omega(\ppr)=
\sum_{i=1}^{n+1}\int_{a^+_i}\Wp(\ppr)\int_{b^+_i}\omega(\ppr)-\int_{b^{+}_i}\Wp(\ppr)\int_{a^+_i}\omega(\ppr),$$
so \be\int_{b^{+}_{n+1}}\Wp(\ppr)=-\sum_{\rm residues}g_{\ppr}\omega(\ppr).\label{eq:reci}\ee
Writing
\be\omega(\ppr)=\frac{\ka(\ppr)\prod_{j=1}^{n}(z-\beta_{j}(\ppr))dz}{w_{+}(\ppr)},\label{eq:omega}\ee
(\ref{eq:reci}) gives that
\be\int_{b^{+}_{n+1}}\Wp(\ppr)=4\ka(\ppr). \label{eq:intbn+1}\ee
For $i=1,\ldots,n+1$,
$$\int_{a^+_i}\dot\omega=\left.\pnu\right|_{\nu=0}\int_{a^+_i}\omega(\ppr)=0,$$
and since $\dot\omega$ is holomorphic, this implies that
$\dot\omega=0$.
Now $$\dot\omega=\left(\dot \ka-\ka(\po)\sum_{j=1}^{n}\frac{\dot
\beta_j}{z-\beta_j(\po)}\right)\frac{\prod_{j=1}^{n}(z-\beta^+_j(\po))dz}{w_{+}(\po)},$$
so this gives $\dot \ka=0$, (and $ \dot\beta_j=0$, $j=1,\ldots,n$) and hence by (\ref{eq:intbn+1}), we have that
$$\left.\pnu\right|_{\nu=0}\int_{b^{+}_{n+1}}\Wp=0.$$

\noindent{\bf\setcounter{enumi}{6}({\it\arabic{enumi}}):} Using the fact that $\dot\Omega_+(\po)=0$, a calculation gives that
\begin{eqnarray}
\ddot{\Omega}_+(\po)%
&=&-\left(\sum_{j=1}^{n}\frac{\ddot{\zeta}^+_j}{z-\zeta^+_j}-\frac{\ddot{\zeta}^+_{n+1}}{z-\mu}-\frac{1}{(z-\m)^2}\right)\Wp(\po).\label{eq:Wddot} \end{eqnarray}

Since
$$\int_{a^+_{n+1}}\ddot\Omega_+(\po)=\left.\ppnu\right|_{\n=0}\int_{a^+_{n+1}}\Omega_+(\po)=0,$$
$\ddot\Omega_+(\po)$ has zero residue at $z=\mu$.

Writing $\Wp(\po)=k(z)dz$, then from equation~(\ref{eq:Wddot}), we have that
$$\ddot\zeta^+_{n+1}k(\mu)=-\frac{dk}{dz}(\mu)$$
and hence using equations~(\ref{eq:Dplusdefn}) and (\ref{eq:Dplus}),%
\begin{equation}\ddot\zeta^+_{n+1}=\frac{1}{2\mu}+\frac{D_+(p)}{\mu^2}+O\left(\frac{1}{\mu^3}\right). \label{eq:ddotzeta}\end{equation}

We can represent the homology classes $a^+_i,b^+_i$, $i=1,\ldots,n$ by loops whose projections to $\cp$ lie in a fixed compact region $K$ that is independent of $\m$. Then for $\mu$ sufficiently large and $z\in K$, we have
\begin{eqnarray}
\frac{1}{z-\mu}&=&-\frac{1}{\mu}-\frac{z}{\mu^2}+O\left(\frac{1}{\mu^3}\right)
\label{eq:zminusm}\\
\frac{1}{(z-\mu)^2}&=&\frac{1}{\mu^2}+\frac{z}{\mu^2}+O\left(\frac{1}{\mu^4}\right)\label{eq:zminusmsq}
\end{eqnarray}

and substituting these, together with (\ref{eq:ddotzeta}), into (\ref{eq:Wddot}), we obtain that in $K$,
\be\ddot\Omega_+(\po)=\left(-\sum_{j=1}^n\frac{\ddot\zeta^+_j}{z-\zeta^+_j}-\frac{1}{2\mu^2}+\frac{2D_+(p)-3z}{2\mu^3}+O\left(\frac{1}{\mu^4}\right)\right)\Wp(\po)\label{eq:ddotW}.\ee

We will thus obtain an asymptotic expression for $\ddot\Omega_+(\po)$ by calculating one for $\ddot\zeta^+_j(\po)$.
To do this, we note that both $\Wp(\po)$ and $\ddot\Omega_+(\po)$ have trivial $a$-periods, so
$$ \frac{3}{2}\int_{a_i(\po)}z\Wp(\po) +\ddot\zeta^+_j\mu^3\sum_{j=1}^n\int_{a^+_i(\po)} \frac{\Omega_+(\po)}{z-\zeta^+_j}=O\left(\frac{1}{\mu}\right) \mbox{, for $i=1,\ldots,n$, }$$
or equivalently
$$ \frac{3}{2}\int_{a_i(p)}z\Wp(p) +\ddot\zeta^+_j\mu^3\sum_{j=1}^n\int_{a^+_i(p)} \frac{\Omega_+(p)}{z-\zeta^+_j}=O\left(\frac{1}{\mu}\right) \mbox{, for $i=1,\ldots,n$. }$$
Recalling the definition of $\Wph$ from equations~(\ref{eq:cj}) and (\ref{eq:Wph}), notice that \be\ddot\zeta^+_j=\frac{c^+_j(p)}{\mu^3}+O\left(\frac{1}{\mu^4}\right).
\label{eq:zetaj}\ee
From (\ref{eq:ddotW}), (\ref{eq:cj}), (\ref{eq:Wph}) and (\ref{eq:zetaj}),
$$\Psi^*_+\ddot\Omega_+= -\frac{1}{2\mu^2}\Wp+\frac{1}{\mu^3}\left(D_+(p)\Wp + \Wph\right)+O\left(\frac{1}{\mu^4}\right),$$
proving all but the last component of \setcounter{enumi}{6}({\it\arabic{enumi}}).

To compute $\left.\ppnu\right|_{\nu=0}\int_{b^+_{n+1}}\Wp$, we again
use reciprocity with $\omega(\po)$. Differentiating (\ref{eq:omega}), and using $\dot\omega=0$ one obtains
\be\ddot\omega=\left(\frac{\ddot\ka}{\ka(\po)}-\frac{1}{(z-\mu)^2}-\sum_{j=1}^{n}\frac{\ddot\beta_j}{z-\beta_j(\po)}\right)\omega(\po).\label{eq:omegaddot}\ee

Taking residues at $z=\mu$, and utilising $${\rm res}_{z=\mu}\omega(\po)=1,\;{\rm res}_{z=\mu}\ddot\omega=0,$$
gives
\be\frac{\ddot\ka}{\ka(\po)}-\sum_{j=1}^{n}\frac{\ddot\beta_j}{z-\beta_j(\po)}={\rm res}_{z=\mu}\frac{\omega(\po)}{(z-\mu)^2}.\label{eq:ddotka}\ee

Thus in order to obtain an asymptotic expression for $\ddot\ka$, we first obtain expressions for $\ka(\po)$, $\beta_j(\po)$ and $\ddot\beta_j$. We make the assumption throughout that $z\in K$. We shall abuse notation and write $\omega(\po)$ for $\Psi_+\!^*(\omega(\po))$, where
\begin{center}$\begin{array}{rccc}
\Psi_+\!^*:&C_+(p)&\longrightarrow&C_+(\po)\\
&(z,w_+(p))&\longmapsto&\left(z,(z-\mu)w_+(p)\right)
\end{array}$
\end{center} is the normalisation map. From (\ref{eq:omega}) and (\ref{eq:zminusm}),
\be\w(\po)=\ka(\po)\left(\frac{-1}{\mu}+O\left(\frac{1}{\mu^2}\right)\right)\frac{\prod_{j=1}^n(z-\beta_j(\po))dz}{w_+(p)},\label{eq:omegaord}\ee
so
$$\int_{a^+_i}\frac{\prod_{j=1}^n(z-\beta_j(\po))dz}{w_+(p)}=O\left(\frac{1}{\mu}\right).$$ $\frac{\prod_{j=1}^n(z-\beta_j(\po))dz}{w_+(p)}$ is moreover a differential of the second kind on $C_+(p)$ whose only singularity is a double pole at $z=\infty$, and it approaches $\frac{z^n}{w_+(p)}$ as $z\ra\infty$. Hence
$$\frac{\prod_{j=1}^n(z-\beta_j(\po))dz}{w_+(p)}+O\left(\frac{1}{\mu}\right)=\Wp(p),$$
from which we conclude that
\be\beta_j(\po)=\zeta_j+\Oone,\,j=1,\ldots,n.\label{eq:betaj}\ee

The fact that res$_{z=\mu}\,\w(\po)=1$ yields
$$\frac{\ka(\po)\prod_{j=1}^n(\mu-\beta_j(\po))}{(w_+(p))(\mu)}=1,$$
and (\ref{eq:betaj}) gives
$$\frac{\prod_{j=1}^n(\mu-\beta_j(\po))}{(w_+(p))(\mu)}=\frac{\prod_{j=1}^n(\mu-\zeta_j)}{(w_+(p))(\m)}+O\left(\frac{1}{\mu^{5/2}}\right)$$
so using (\ref{eq:Dplus}),
$$\ka(\po)\left(\mu^{-1/2}+D_+(p)\mu^{-3/2}+O(\mu^{-5/2})\right)=1$$
and therefore
$$\ka(\po)=\mu^{1/2}-D_+(p)\mu^{-1/2}+O(\mu^{-3/2}).$$
This, along with (\ref{eq:omegaord}) gives that
\be\omega(\po)=\mu^{-1/2}\Omega_+(p)+O(\mu^{-3/2}).\label{eq:omega2}\ee
Since $\int_{a^+_i}\omega(\po)=0$ and $\int_{a^+_i}\ddot\omega=0$, (\ref{eq:omegaddot}), (\ref{eq:omega2}) and (\ref{eq:zminusmsq}) yield, for $i=1,\ldots,n$,
\begin{eqnarray*}
\sum_{j=1}^{n}\ddot\beta_j\int_{a^+_i}\frac{\omega(\po)}{z-\beta_j(\po)}&=&\int_{a^+_i}
\frac{\omega(\po)}{(z-\mu)^2}\\
&=&\frac{2}{\mu^{7/2}}\int_{a^+_i}z\Omega_+(\po) + O\left(\frac{1}{\mu^{9/2}}\right).
\end{eqnarray*}
For each $i$, $\int_{a^+_i}z\Omega_+(\po)$ is independent of $\mu$ and so
\be \sum_{j=1}^{n}\ddot\beta_j\mu^{1/2}\int_{a^+_i}\frac{\omega(\po)}{z-\beta_j(\po)}=O\left(\frac{1}{\mu^3}\right). \label{eq:order3}\ee
But by (\ref{eq:betaj}) and (\ref{eq:omega2}),
$$\left(\mu^{1/2}\int_{a^+_i}\frac{\omega(\po)}{z-\beta_j(\po)}\right)^i_j=
\left(\int_{a^+_i}\frac{\Omega_+(p)}{z-\zeta_j}\right)^i_j+O\left(\frac{1}{\mu}\right),$$
and the matrix on the right is invertible and has no dependence on $\mu$.
(\ref{eq:order3}) hence implies that for $j=1,\ldots,n$,
$$\ddot\beta_j(\po) = O\left(\frac{1}{\mu^3}\right).$$

We now have the asymptotics for $\ka(\po)$, $\beta_j(\po)$ and $\ddot\beta_j$ that we desired earlier, and substituting them into
\begin{eqnarray*}
\frac{\ddot\ka}{\ka(\po)}-\sum_{j=1}^{n}\frac{\ddot\beta_j}{z-\beta_j(\po)}&=&
{\rm res}_{z=\mu}\frac{\omega(\po)}{(z-\mu)^2}\hspace{35 mm}\mbox{ from (\ref{eq:ddotka})}\\
&=&\left.\frac{1}{2}\frac{d^2}{dz^2}\right|_{z=\mu} \left(\frac{\ka(\po)\prod_{j=1}^{n}(z-\beta_j(\po))}{w_+(p)}\right)
\end{eqnarray*}
gives
$$\ddot\ka=\frac{3}{8\mu^{3/2}}+\frac{9D_+(p)}{8\mu^{5/2}}+O\left(\frac{1}{\mu^{7/2}}\right)$$
which upon substitution into (\ref{eq:intbn+1}) completes the proof of \setcounter{enumi}{6}({\it\arabic{enumi}}) . \begin{flushright}$\Box$\end{flushright}

We now have asymptotic expressions for each row of the $2n+5\times 2n+5$ matrix $H(\m)$ in (\ref{eq:h}), which we wish to show is non-singular in the limit as $\m\ra\infty$ along $L_\epsilon$ of figure~\ref{fig:hdomain}. The inductive assumption and Lemma~\ref{thm:asym} tell us that columns \mbox{$1,\ldots,n+2,n+4,\ldots,2n+4$} of the the first $2n+3$ rows of $H(\m)$ are linearly independent, and that its $2n+4^{th}$ row $(\pmu I_+(\po);\pmu I_-(\po))$ is

\vspace{4mm}
\noindent${\ds(\underbrace{0,\ldots,0}_{\mbox{\scriptsize $n+2$ zeros}},2\im(\m^{-1/2}+ D_{+}(p)\m^{-3/2})+O(\m^{-5/2});}$
\vspace{-6mm}\begin{flushright}${\ds\underbrace{0,\ldots,0}_{\mbox{\scriptsize $n+1$ zeros}},2(\m^{-1/2}+D_{-}(p)\m^{-3/2})+O(\m^{-5/2})).}$\end{flushright}
Note that the two non-zero entries in this row have leading terms differing only by multiplication by $\im$. We find a linear combination of the rows of $H(\m)$ that equals
${\displaystyle(\underbrace{0,\ldots,0}_{\mbox{\scriptsize $n+2$ zeros}},
\im\m^{-5/2}(4\eta^+(p)-5 D_+(p))+O(\m^{-7/2});}$\vspace{-6mm}\begin{flushright}$\displaystyle{\underbrace{0,\ldots,0}_{\mbox{\scriptsize
$n+1$ zeros}},
\m^{-5/2}(4\eta^-(p)-5D_-(p))+O(\m^{-7/2}))},$\end{flushright} where $\eta^\pm(p)$ are defined in (\ref{eqn:etas}).

Our matrix is non-singular if
$$\lim_{\m\ra\infty}  4\eta^+(p)-5D_+(p) \neq \lim_{\m\ra\infty}
4\eta^-(p)-5D_-(p)$$
where the limits are taken along $L_\epsilon$. This is the ``extra condition'' referred to earlier, and we will modify the statement we prove  by  induction to ensure that it is satisfied. First, we find the linear combination yielding this condition.

From  Lemma~\ref{thm:asym} we have that
\be\hspace*{-103mm}\left(\ppnu I_+(\po);\ppnu I_-(\po)\right)=\label{eq:ppnu}\ee
$\displaystyle{\left(\left(\frac{D_+(p)}{\m^3}-\frac{1}{2\m^2}\right)I_+(p)
-\frac{1}{\m^3}\left(\widehat I_+(p)+O\left(\frac{1}{\m}\right)\right),
\frac{3\im}{2\m^{3/2}}+\frac{9\im D_{+}}{2\m^{5/2}}+O\left(\frac{1}{\m^{7/2}}\right);\right.}$    \begin{flushright}
${\displaystyle\left.\left(\frac{D_-(p)}{\m^3}-\frac{1}{2\m^2}\right)I_-(p)
-\frac{1}{\m^3}\left(\widehat I_-(p)+O\left(\frac{1}{\m}\right)\right),
\frac{3}{2\m^{-3/2}}+\frac{9 D_{-}}{2\m^{5/2}}+O\left(\frac{1}{\m^{7/2}}\right)\right)}$
\end{flushright}

By the induction hypothesis, there are unique \mbox{$\eta^{\pm}(p)$},
\mbox{$\chi(p)$}, \mbox{$\xi_i(p)$}, \mbox{$i=1,\ldots2n$} such that
\begin{eqnarray}\hspace*{-10mm}(\widehat{I}_+(p);\widehat I_-(p)) & = & \eta^+(p)(I_+(p);0)+\eta^-(p)(0;I_-(p))+\chi(p)\frac{\partial}{\partial
R}({I}_+(p);I_-(p))\\
&&+\sum_{i=1}^{2n}\xi_{i}(p)\pli(I_+(p);I_-(p)).\label{eqn:etas}
\end{eqnarray}
Thus by (\ref{eq:ppnu}) and Lemma~\ref{thm:asym} there are \mbox{$\tilde\eta^{\pm}(p)=\eta^{\pm}(p)+O(\m^{-1})$}, \mbox{$\tilde\chi(p)=\chi(p)+O(\m^{-1})$},
\mbox{$\tilde\xi_i(p)=\xi_i(p)+O(\m^{-1})$} such that

\vspace{4mm}
\noindent${\ds(\ppnu I_+(\po);\ppnu I_-(\po))
+(\frac{1}{2\m^2}+\frac{\tilde\eta^+(p)-D_+(p)}{\m^3})\left(I_+(\po);0\right)} + $\\*
${\ds(\frac{1}{2\m^2}+\frac{\tilde\eta_-(p)-D_+(p)}{\m^3})\left(0;I_-(\po)\right)+\frac{\tilde\chi(p)}{\m^3}\frac{\partial}{\partial R}(I_+(p);I_-(p))
+\sum_{i=1}^{2n}\frac{\tilde\xi_{i}(p)}{\m^3}\pli(I_+(p);I_-(p))}$\\*
\vspace{-4mm}
\be\hspace*{40mm}=(\underbrace{0,\ldots,0}_{\mbox{\scriptsize $n+2$ zeros}},g_+(\po);
\underbrace{0,\ldots,0}_{\mbox{\scriptsize $n+1$ zeros}},g_-(\po)) \label{eq:lincomb}\ee
where
\noindent\begin{eqnarray*}
\frac{g_+(\m)}{\im}&\hspace*{-6pt}=&\hspace{-6pt}\left(\frac{3}{2\m^{3/2}}+\frac{9 D_{+}(p)}{2\m^{5/2}}
+(\frac{1}{2\m^2}+\frac{\eta^+(p) - D_+(p)}{\m^3})(4\m^{1/2}
-\frac{4 D_{+}(p)}{\m^{1/2}})\right)+O\left(\frac{1}{\m^{7/2}}\right)\\
&\hspace*{-6pt}=&\hspace*{-6pt}\frac{7}{2\m^{3/2}}-\frac{3 D_+(p)}{2\m^{5/2}}+\frac{4\eta^+(p)}{\m^{5/2}}
+O\left(\frac{1}{\m^{7/2}}\right),
\end{eqnarray*}
and similarly
$$g_-(\m)=\frac{7}{2\m^{3/2}}-\frac{3 D_-(p)}{2\m^{5/2}}+\frac{4\eta^-}{\m^{5/2}}
+O\left(\frac{1}{\m^{7/2}}\right),$$
so adding the appropriate multiple of $\pmu(I_+(\po),I_-(\po))$ to (\ref{eq:lincomb}), we have
\noindent\begin{eqnarray}
l(\po)&:=&\left(\ppnu I_+(\po);\ppnu I_-(\po)\right)+
\left(\frac{1}{2\m^2}+\frac{\tilde\eta^+(p)-D_+(p)}{\m^3}\right)\left(I_+(\po);0\right)
\nonumber\\
&&+(\frac{1}{2\m^2}+\frac{\tilde\eta_-(p)-D_-(p)}{\m^3})\left(0;I_-(\po)\right)
+\frac{\tilde\xi}{\m^3}\frac{\partial}{\partial R}(I_+(p);I_-(p))\nonumber\\
&&+\frac{1}{\m^3}\sum_{i=1}^{2n}\tilde\xi_{i}(p)\pli(I_+(p);I_-(p))
-\frac{7}{4\m}\pmu(I_+(\po),I_-(\po))\nonumber\\
&=&(\underbrace{0,\ldots,0}_{\mbox{\scriptsize $n+2$ zeros}},
\frac{\im(4\eta^+(p)-5 D_+(p))}{\m^{5/2}}+O\left(\frac{1}{\m^{7/2}}\right);\nonumber\\
&&\hspace*{30mm}\underbrace{0,\ldots,0}_{\mbox{\scriptsize $n+1$ zeros}},
\frac{(4\eta^-(p)-5D_-(p))}{\m^{5/2}}+O\left(\frac{1}{\m^{7/2}}\right)\label{eq:lincomb2}.\end{eqnarray}

We are led therefore, to modify the statement we prove by induction to include the assumption that

$$ \lim_{\m\ra\infty}  4\eta^+(p)-5D_+(p) \neq  \lim_{\m\ra\infty}
4\eta^-(p)-5D_-(p),$$%
where the limits are taken along $L_\epsilon$. Of course this modification needs to be such that it is preserved under the induction step. With this in mind, we define $\eta^\pm(\po)$ by the condition that\\
\vspace{5mm}
${\ds (\widehat{I}_+(\po);\widehat{I}_-(\po)) - \eta^+(\po)(I_+(\po);0) +
\eta^-(\po)(0;I_-(\po))\in}$
\vspace{-15mm}\begin{flushright} ${\ds \mbox{span}\left\{\pR(I_+(\po);I_-(\po)),
\pli (I_+(\po);I_-(\po)),\pmu (I_+(\po);I_-(\po)), \ppnu (I_+(\po);I_-(\po))\right\}},$\end{flushright}
\vspace{2mm}
and calculate the relationship between $\eta^+(\po)-\eta^-(\po)$ and
$\eta^+(p)-\eta^-(p)$.

\begin{lemma}
As $\m\ra\infty$ along $L_\epsilon$,
\begin{eqnarray*}
\widehat I_+(\po)&=&(\widehat I_+(p),2\im\m^{3/2}+
(6\im D_+(p)-4\im\eta^+(p))\m^{1/2} + O(\m^{-1/2})),\\
\widehat I_-(\po)&=&(\widehat I_-(\po),2\m^{3/2}+(6D_-(p)-4\eta_-(p))\m^{1/2} + O(\m^{-1/2})).
\end{eqnarray*}\label{thm:widehatIpm}\end{lemma}
\noindent {\bf Proof of Lemma~\ref{thm:widehatIpm}:}
For $\m$ sufficiently large, we may assume that $\m\neq\zeta^\pm_j$,\\$j=1,\ldots,n$. Then
$$\widehat\Omega_+(\po):=\frac{3}{2}z\Wp(\po)
+\sum_{j=1}^{n+1}\frac{c^+_j(\po) \Wp(\po)}{z-\zeta^+_j},$$
where the $c^+_j(\po)$, $j=1,\ldots,n+1$ are determined by the (non-singular)
system of equations
$$\frac{3}{2}\int_{a^\pm_i}z\Wpm(\po)+\sum_{j=1}^{n+1}c^\pm_j(\po)\int_{a^\pm_i}\frac{\Wpm(\po)}{z-\zeta^\pm_j}=0,\:i,\,j=1,\ldots,n+1.$$
Taking $i=n+1$ we quickly see that \mbox{$c^+_{n+1}(\po)=0$} and \mbox{$c^+_j(\po)=c^+_j(p)$,} for \mbox{$j=1,\ldots,n$}, so, not surprisingly,
$$\Psi_+\!^*(\widehat \Omega_+(\po))=\widehat \Omega_+(p).$$
The argument for $\widehat\Omega_-(\po)$ is similar, and
together they prove all but the last components of the lemma above. For these, we again let  $\Gamma$ denote the circle $|z|=\m$, traversed clockwise.
For $\m$ sufficiently large,
\begin{eqnarray*}\int_{b^+_{n+1}}\widehat\Omega_+(\po)
&=&-\int_{\Gamma}\left(\frac{3}{2}z+\sum_{j=1}^{n}\frac{c^+_j(p)}{z-\zeta^+_j}\right)
\frac{\prod_{k=1}^n(z-\zeta^+_k)dz}{\sqrt{(z-R)\prod_{i=1}^{2n} (z-\lambda_i)}}\\
&=&2\im\m^{3/2}+(6\im D_+(p)-4\im\eta^+(p))\m^{1/2}+O(\m^{-1/2}),
\end{eqnarray*}
and similarly for \mbox{$\int_{b^-_{n+1}(\po)}\widehat\Omega_-(\po)$}.
\begin{flushright}$\Box$
\end{flushright}

From Lemma~\ref{thm:asym},

\noindent${\ds (\widehat I_+(\po);\widehat I_-(\po))-\eta^+(p)(I_+(\po);0)
-\eta^-(p)(0;I_-(\po))-\chi(p)\pR(I_+(\po);I_-(\po))}$\\*
${\ds -\sum_{i=1}^{2n}\xi_i(p)\pli(I_+(\po);I_-(\po))}$\\*
\hspace*{5mm}${\ds=(0,2\im\m^{3/2}+(6\im D_+(p)-4\im\eta^+(p))\m^{1/2}+O(\m^{-1/2});}$
\vspace{-5mm}\begin{flushright}${\ds 0, 2\m^{3/2}+(6D_-(p)-4\eta^-(p))\m^{1/2}+O(\m^{-1/2})}$\end{flushright}
\vspace{-5mm}
\be\hspace*{-77mm}=\Lambda \pmu(I_+(\po);I_-(\po)) + \Upsilon l(\po),\label{eq:LamUps}\ee
where $l(\po)$ is defined in (\ref{eq:lincomb2}), and by
Lemma~\ref{thm:widehatIpm}, we have that $\Lambda$ and $\Upsilon$ are defined by the equations
\begin{flushleft}${\ds \left(\begin{array}{cc}
\frac{2\im}{\m^{1/2}}(1 + \frac{D_+(p)}{\m})+O(\frac{1}{\m^{5/2}}) &
\frac{\im(4\eta^+(p)-5D_+(p))}{\m^{5/2}}+O(\frac{1}{\m^{7/2}})\\
\frac{2}{\m^{1/2}}(1 + \frac{D_-(p)}{\m})+O(\frac{1}{\m^{5/2}}) &
\frac{4\eta^-(p)-5D_-(p)}{\m^{5/2}}+O(\frac{1}{\m^{7/2}})
\end{array}\right)
\left(\begin{array}{c}
\Lambda\\
\Upsilon
\end{array}\right)}$\end{flushleft}\nopagebreak
\begin{flushright}${\ds =\left(\begin{array}{c}
2\im\m^{3/2}+ (6\im D_+(p)-4\im\eta^+(p))\m^{1/2} + O(\frac{1}{\m^{1/2}})\\
2\m^{3/2}+ (6 D_-(p)-4\eta^-(p))\m^{1/2} + O(\frac{1}{\m^{1/2}})
\end{array}\right)}$\end{flushright}
so
\begin{eqnarray*}
\Lambda&=&\m^2+O(\m),\\
\Upsilon&=&\frac{4(\eta^+(p)-\eta^-(p))-4(D_+(p)-D_-(p))}{5(D_+(p)-D_-(p))-4(\eta^+(p)-\eta^-(p))}\m^3+O(\m^2).\end{eqnarray*}
Using (\ref{eq:lincomb2}) and (\ref{eq:LamUps}) then,
$$\eta^+(\po)=\eta^+(p)+\Upsilon\left(\frac{1}{2\m^2}+\frac{\eta^+(p)-D_+(p)}{\m^3}+O\left(\frac{1}{\m^4}\right)\right),$$
$$\eta^-(\po)=\eta^-(p)+\Upsilon\left(\frac{1}{2\m^2}+\frac{\eta^-(p)-D_-(p)}{\m^3}+O\left(\frac{1}{\m^4}\right)\right),$$
so
\begin{eqnarray}
\lefteqn{\eta^+(\po)-\eta^-(\po)}\nonumber\\
&=&\eta^+(p)-\eta^-(p)\nonumber\\
&&\hspace*{-7mm}+\frac{4(\eta^+(p)-\eta^-(p))-4(D_+(p)-D_-(p))}{5(D_+(p)-D_-(p))-4(\eta^+(p)-\eta^-(p))}\left(\eta^+(p)-\eta^-(p)-(D_+(p)-D_-(p)\right)
+O\left(\frac{1}{\m}\right)\nonumber\\
&=&\frac{(D_+(p)-D_-(p))\left(3(\eta^+(p)-\eta_-(p))-4(D_+(p)-D_-(p))\right)}
{4(\eta^+(p)-\eta_-(p))-5(D_+(p)-D_-(p))}.\label{eq:etapo}
\end{eqnarray}
Defining $T_p$ to be the linear fractional transformation
$$T_p:x\mapsto (D_+(p)-D_-(p))\frac{3x-4(D_+(p)-D_-(p))}{4x-5(D_+(p)-D_-(p))},$$
(\ref{eq:etapo}) is the statement that
$$T_p(\eta^+(p)-\eta^-(p))=\eta^+(\po)-\eta^-(\po).$$
Moreover, we know that $D_\pm(\po)=D_\pm(p)$ and hence
$$T_{\po}=T_p.$$
The statement then, that we prove by induction, is:
\begin{thm}
    For each positive integer $m$ and integer $n$ with $0\leq n\leq m$ there exists $p\in M_{n,\R}$ such that
\begin{enumerate}
\item{$\zeta^+_j(p)$, $j=1,\ldots,n$ are pairwise distinct, as are
$\zeta^-_j(p)$, $j=0\ldots n$,}
\item $\mathbb R^{2n+3}$ is spanned by the vectors $(I_{+}(p),0)$, $(0,I_{-}(p))$,
$\frac{\partial}{\partial R}\left(I_{+}(p), I_{-}(p)\right)$
 and $\pli\left(I_{+}(p), I_{-}(p)\right)$, $i=1,\ldots,2n$,
\item{$5\left(D_{+}(p)-D_{-}(p)\right)+4T^{k}_p\left(\eta_{+}(p)-\eta_{-}(p)\right)\neq
0$, for $0\leq k\leq m-n$.}
\end{enumerate}
\label{thm:oddextra}
\end{thm}%
For the convenience of the reader, we reiterate here the definitions of many of the objects appearing in Theorem~\ref{thm:oddextra}, so that it  may be read immediately after the statement of Theorem~\ref{thm:odd}, without recourse to the above arguments.
Assume that \setcounter{enumi}{1}({\it\arabic{enumi}}) and \setcounter{enumi}{2}({\it\arabic{enumi}}) of Theorem~\ref{thm:oddextra} hold for
$p\in M_{n,\R}$. Then the differentials
$\frac{\Wpm(p)}{z-\zeta^\pm_j}$ are a basis for the holomorphic
differentials on $\Cpm(p)$. Thus we may define
$c^\pm_j(p)$ by the equations
$$\frac{3}{2}\int_{a^\pm_i}z\Wpm(p)+\sum_{j=1}^{n}c^\pm_j(p)\int_{a^\pm_i}\frac{\Wpm(p)}{z-\zeta^\pm_j}=0,\:i,j=1,\ldots,n.$$
Let
$$\Wpmh(p):=\frac{3}{2}z\Wpm(p)+\sum_{j=1}^{n}c^\pm_j(p)\frac{\Wpm(p)}{z-\zeta^\pm_j},$$
 $$\widehat{I}_+(p):=\im\left(\int_{c_1}\widehat{\Omega}_+(p),\int_{c_{-1}}\widehat{\Omega}_+(p),\int_{b^+_1}\widehat{\Omega}_+(p),\ldots,\int_{b^+_n}\widehat{\Omega}_+(p)\right)$$
and
$$\widehat{I}_-(p):=\left(\int_{b^-_0}\widehat{\Omega}_-(p),\int_{b^-_1}\widehat{\Omega}_-(p),\ldots,\int_{b^-_n}\widehat{\Omega}_-(p)\right).$$
Define
\mbox{$\eta^{\pm}(p)$}, \mbox{$\xi(p)$}, \mbox{$\xi^\pm_i(p)$}, \mbox{$i=1,\ldots2n$} by
\begin{eqnarray*}
(\widehat{I}_+(p),\widehat I_-(p))&=&
\eta^+(p)(I_+(p),0)+\eta^-(p)(0,I_-(p))+\xi(p)\frac{\partial}{\partial
R}({I}_+(p),I_-(p))\\
&&+\sum_{i=1}^{2n}\xi_{i}(p)\pli(I_+(p),I_-(p)).\end{eqnarray*}

Put
$$D_{+}(p):=\frac{1}{2}\left(R+\sum_{i=1}^{2n}\lambda_{i}\right)-\sum_{j=1}^n\zeta^+_j,$$
$$D_{-}(p):=\frac{1}{2}\left(R+\sum_{i=1}^{2n}\lambda_{i}\right)-\sum_{j=0}^n\zeta^-_j$$
and let $T_p$ be the linear fractional transformation
$$T_p:x\mapsto (D_+(p)-D_-(p))\frac{3x-4(D_+(p)-D_-(p))}{4x-5(D_+(p)-D_-(p))}.$$

\noindent{\bf Proof of Theorem~\ref{thm:oddextra}:} Fix $m$, and for $n<m$ suppose $p\in M_{n,\R}$ satisfies the conditions of Theorem~\ref{thm:oddextra}. By the above arguments, the set of $\m\in(-2,2)$ such that
\begin{enumerate}
\item[\setcounter{enumi}{1}({\roman{enumi}})] for all $\epsilon\in (0,\min_{i = 1,\ldots,n}|\lambda_i+2|)$, \mbox{$ h_\epsilon(\m)\neq 0$},

\item[\setcounter{enumi}{2}({\roman{enumi}})] for $j=1,\ldots n$, $\m\neq\zeta^+_j$ and

\item[ \setcounter{enumi}{3}({\roman{enumi}})]for $j=0,\ldots n$, $\m\neq\zeta^-_j$
\end{enumerate}
is dense  in $(-2,2)$.

Take such a $\m$. Then $\po=(p,\m,0)$ satisfies
Theorem~\ref{thm:oddextra}, where in \setcounter{enumi}{2}({\it\arabic{enumi}}) we replace \mbox{$\frac{\partial}{\partial\lambda_{2n+2}}(I_+(\po); I_-(\po))$} by \mbox{$\ppnu(I_+(\po); I_-(\po))$}.
Then for $\n$ small, utilising (\ref{eq:nusmall}),
$$\eta^\pm(p,\m,\n) = \eta^\pm(\po) + O(\n),$$
$$D_\pm(p,\m,\n) = D_\pm(\po) +O(\n)$$
and
$$T_{(p,\m,\n)}=T_{\po}+O(\n),$$
we conclude that $(p,\m,\n)$ satisfies Theorem~\ref{thm:oddextra}. It remains to show the existence of $p\in M_{0,\R}$ verifying \setcounter{enumi}{1}({\it\arabic{enumi}}) and \setcounter{enumi}{2}({\it\arabic{enumi}}) of Theorem~\ref{thm:oddextra}, and such that for no $k\geq 0$ do we have $5\left(D_{+}(p)-D_{-}(p)\right)+4T^{k}_p\left(\eta_{+}(p)-\eta_{-}(p)\right)=0$.

\subsection{Genus One ($n=0$)}
\begin{figure}[h]
\hspace*{1.5cm}\includegraphics{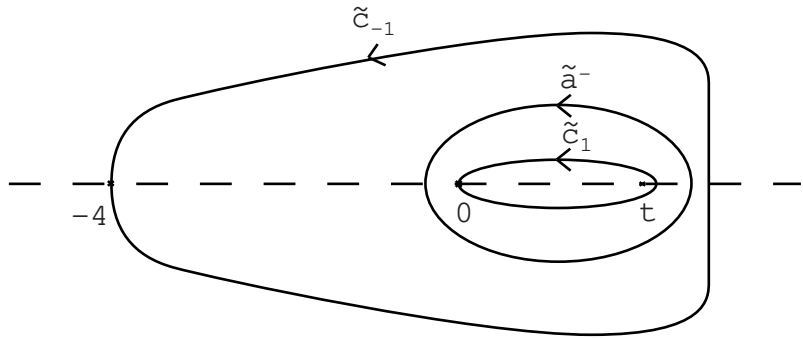}
\caption{ The curves $\tilde a^-$ and $\tilde c_{\pm 1}$, for $n=0$.\label{fig:oneCpm}}
\end{figure}

We consider pairs $C_+=C_+(R)$ and $C_-=C_-(R)$ given by
$$w_+^2=(z-R)$$
and
$$w_-^2=(z+2)(z-2)(z-R)$$ respectively, where $R>2$. Writing 
$\pi_\pm:(z,w_\pm)\mapsto z$ for the projections to $\cp$, the fibre product 
of these is the genus one curve $X=X(R)$, given by
$$y^2=x(x-r)(x-\frac{1}{r})\mbox{, where } r+\frac{1}{r}=R.$$
\begin{lemma}
There exists a $p\in M_{0,\R}$ such that
\begin{enumerate}
\item $\mathbb{R}^{3}$ is spanned by the vectors $(I_+(p),0),\,(0,I_-(p))$ and $\pR(I_+(p),I_-(p))$,
\item for all $k\geq 0$, $5(D_+(p)-D_-(p))+4T^k(p)(\eta^+(p)-\eta^-(p))\neq 0$.
\end{enumerate}
\label{lem:one}
\end{lemma}
\smallskip\noindent{\bf Proof:} The natural limit to consider is $r\ra 1$, i.e. $R=r+1/r \ra 2$, which suggests setting $\zeta:=z+2$, $t:=R-2$. Then $C_+(t)$ is given by $$w_+^2:=\z-t,$$
and $C_-(t)$ by
$$w_-^2:=\z(\z-t)(\z+4).$$
For each $t>0$, choose $c_1(t), c_{-1}(t)$ and $a^-(t)$ as shown in figure~\ref{fig:oneCpm}.
We write $\Wm(t)=\frac{(\z-s(t))d\z}{w_-}$, where $s(t)$ is defined by the condition $\int_{a^-(t)}\Wm(t) = 0$. Using 
$$\int_{a^-(t)}\Wm(t) = 0\;\mbox{ and }\; \frac{\partial}{\partial t}\int_{a^-(t)}\Wm(t) = 0, $$
one obtains 
\be s(t)=t+O(t^2).\label{eq:s}\ee 
Thus
\begin{eqnarray}
I_-(t)&=&\int_{b_-(t)}\Wm(t)\\
 &=&2\int_{-4}^0\frac{(\z-t)d\z}{\sqrt{\z(\z+4)(\z-t)}}+O(t^2)\\
&=&8+O(t).
\label{eq:oneIminus}
\end{eqnarray}
Now $I_+(t)=\im(\int_{c_1(t)}\Wp(t),\,\int_{c_{-1}(t)}\Wp(t))$, where $c_1(t)$ is a path in $C_+(t)$ joining the two points with $\z=0$, and $c_{-1}(t)$ is one joining the two points with $\z=-4$, both beginning at points with $\frac{w_+}{\im}<0$. Then  
\begin{eqnarray}
I_+(t)&=&2\im \left(\int_0^t\frac{d\z}{\sqrt{\z-t}},\int_{-4}^t\frac{d\z}{\sqrt{\z-t}}\right)\\
&=&(4t^{1/2},4(4+t)^{1/2})
\label{eq:oneIplus}
\end{eqnarray} and
$$\frac{\partial I_+(t)}{\partial t} = (2t^{-1/2},2(4+t)^{-1/2})$$ so we see that condition (1) of Lemma~\ref{lem:one} is satisfied for all $t>0$, in other words for all $R>2$. 

We proceed then to calculate $\lim_{t\ra 0}D_\pm(t)$ and $\lim_{t\ra 0}\eta^\pm(t)$. We have 
\bes
\lim_{t\ra 0}D_+(t)&=&\lim_{t\ra 0}\frac{1}{2}(t+2)\\
&=&1
\ees
and
\bes
\lim_{t\ra 0}D_-(t)&=&\lim_{t\ra 0}\frac{1}{2}(t+2)-(2+s(t))\\
&=&-1.
\ees
Also,
$$\widehat\Omega_+(t)=\frac32z\Wp(t) = \frac{3(\z+2)d\z}{2\sqrt{\z-t}},$$ so
\be
\widehat I_+(t)&=&\im\left(\int_{c_1(t)}\widehat\Omega_+(t),\,\int_{c_{-1}(t)}\widehat\Omega_+(t)\right)\\
&=&\im\left(3\int_0^t\frac{(\z+2)d\z}{\sqrt{\z-t}},\,3\int_{-4}^t\frac{(\z+2)d\z}{\sqrt{\z-t}}\right)\\
&=&(12t^{1/2}+8t^{3/2},\,12(4+t)^{1/2}+8(4+t)^{3/2}).
\label{eq:oneIhatplus}
\ee
Recall that $\widehat\Omega_-(t)=\frac{3}{2}(\z+2)\Wm(t)+\frac{c_-(t)\Wm(t)}{\z-s(t))}$, where $c_-(t)$ is defined by
\be \frac{3}{2}\int_{a^-(t)}\frac{(\z+2)(\z-s(t))d\z}{\sqrt{\z(\z+4)(\z-t)}}+c_-(t)\int_{a^-(t)}\frac{d\z}{\sqrt{\z(\z+4)(\z-t)}}=0.\label{eq:cminus}\ee
Let $\int_{a^-(t)}\frac{d\z}{\sqrt{\z(\z+4)(\z-t)}}=a_0+a_1t+O(t^2)$. Then
\bes
a_0&=&2\pi\im\mbox{Res}_{\z=0}\frac{1}{\sqrt{\z(\z+4)(\z-t)}}\\
&=&\pi\im,
\ees
and
\bes
a_1&=&2\pi\im\mbox{Res}_{\z=0}\left.\frac{\partial }{\partial t }\right|_{t=0}\frac{1}{\sqrt{\z(\z+4)(\z-t)}}\\
&=&\frac{-\pi\im}{16},
\ees
so $$\int_{a^-(t)}\frac{d\z}{\sqrt{\z(\z+4)(\z-t)}}=\pi\im -\frac{\pi\im}{16}t + O(t^2),$$
and similarly
$$\int_{a^-(t)}\frac{\z d\z}{\sqrt{\z(\z+4)(\z-t)}}=\pi\im t +O(t^2),$$
$$\int_{a^-(t)}\frac{\z^2 d\z}{\sqrt{\z(\z+4)(\z-t)}}=O(t^2).$$
Substituting these and (\ref{eq:s}) into (\ref{eq:cminus}), one obtains
$$c_-(t)=O(t^2).$$
From this and (\ref{eq:s}),
\be
\widehat I_-(t)&=&3\int_{-4}^0 \frac{(\z +2)d\z}{\sqrt{\z +4}} +O(t)\\
&=&-8+O(t).
\label{eq:oneIhatminus}
\ee
From (\ref{eq:oneIminus}), (\ref{eq:oneIplus}), (\ref{eq:oneIhatminus}) and (\ref{eq:oneIhatplus}) then
\begin{eqnarray*}
\left(\begin{array}{cc}
I_+(t)&0\\
0&I_-(t)\\
\frac{\partial I_+(t)}{\partial t}&\frac{\partial I_-(t)}{\partial t}\\
\widehat I_+(t)&\widehat I_-(t)
\end{array}\right)&=&\left(\begin{array}{ccc}
4t^{1/2}&4(4+t)^{1//2}&0\\
0&0&8+O(t)\\
2t^{-1/2}&2(4+t)^{-1/2}&O(1)\\
12t^{1/2}+8t^{3/2}&12(4+t)^{1/2}+8(4+t)^{3/2}&-8+O(t)
\end{array}\right).
\end{eqnarray*}
Upon multiplication of its third row by $t$, its first column by $2t^{-1/2}$ and its second column by \(2(4+t)^{-1/2}\) this matrix becomes
\[\left(\begin{array}{ccc}
2&2&0\\
0&0&8\\
1&0&0\\
6&22&-8
\end{array}\right)+O(t).\]
Since \(-11(2,2,0)+1(0,0,8)+16(1,0,0)+1(6,22,-8)=(0,0,0)\),
then recalling that \(\eta^\pm(t)\) are defined by the condition
$$(\widehat I_+(t),\widehat I_-(t))+\eta^+(t)(I_+(t),0)+\eta^-(t)(0,I_-(t))\in {\rm span} \left\{\frac{\partial}{\partial t}(I_+(t),I_-(t))\right\},$$
we conclude that
$$\lim_{t\ra 0}\eta^+(t)=-11$$ and
$$\lim_{t\ra 0}\eta^-(t)=1.$$
The linear fractional transformation $T_t$ is defined by
$$T_t:u\mapsto \frac{-(D_-(t)-D_+(t))(3u+4(D_-(t)-D_+(t)))}{4u+5(D_-(t)-D_+(t))},$$
so letting $T:=\lim_{t\ra 0}T_t$,
$$T:u\ra \frac{3u-8}{2u-5}.$$
This has a  unique fixed point ($u=2$) and so is conjugate to a translation, in fact denoting the map $u\mapsto\frac{1}{u-2}$ by $S$, we have
$$STS^{-1}:u\mapsto u-2.$$
Now
\[ \begin{array}{lccc}
&4T^k(\lim_{t\ra 0}(\eta^-(t)-\eta^+(t)))&=&5(\lim_{t\ra 0}(D_+(t)-D_-(t)))\\
\iff &T^k(12)&=&\frac52\\
\iff &(STS^{-1})^k(\frac{1}{10})&=&2\\
\iff &\frac{1}{10}-2k&=&2,
\end{array}\]
which is clearly false for all integers $k\geq 0$. 
Thus for $t>0$ sufficiently small, the Lemma holds. \begin{flushright}$\Box$ \end{flushright}

\section{Even Genera}
Let $C_\pm=C_\pm(\lambda_1,\ldots \lambda_{2n})$ be the curves given by
$$w_\pm^2=(z\pm 2)\prod_{i=1}^{2n}(z-\lambda_i),$$
where we assume that the sets \mbox{$\{\lambda_1,\lambda_{2}\},\ldots,\{\lambda_{2n-1},\lambda_{2n}\}$} are mutually disjoint and that \mbox{$\lambda_i\neq\pm 2$} for \mbox{$i=1,\ldots,2n$}. Denote by \mbox{$\pi_{\pm}:\Cpm\ra\mathbb{C}P^1$} the respective projections \mbox{$(z,w_\pm)\mapsto z$} to the Riemann sphere, and define a real structure $\rho_\pm$ on $\Cpm$ by 
$$\rho_\pm(z,w_\pm)= (\bar z,\pm\bar w_\pm).$$

Let \mbox{$\tilde a_1\ldots,\tilde a_n, \tilde c_1, \tilde c_{-1}$} be closed curves in  
\mbox{$\cp-\{\lambda_1,\ldots,\lambda_{2n}\}$}, such that 
\begin{enumerate}
\item $\tilde a_i$ has winding number one around $\lambda_i$ and $ \lambda_{n+i}$, winding number zero around the other branch points of $C_\pm$, and does not pass through $2$ or $-2$,
\item $\tilde c_1$ begins and ends at $z=2$, has winding number one around $-2$ and each $\lambda_i$,\\ $i=1,\ldots,2n$,
\item $\tilde c_{-1}$ begins and ends at $z=-2$, has winding number one around $2$ and each $\lambda_i$,\\ $i=1,\ldots,2n$,
\end{enumerate} 

\begin{figure}[h]
\hspace*{1.5cm}\includegraphics{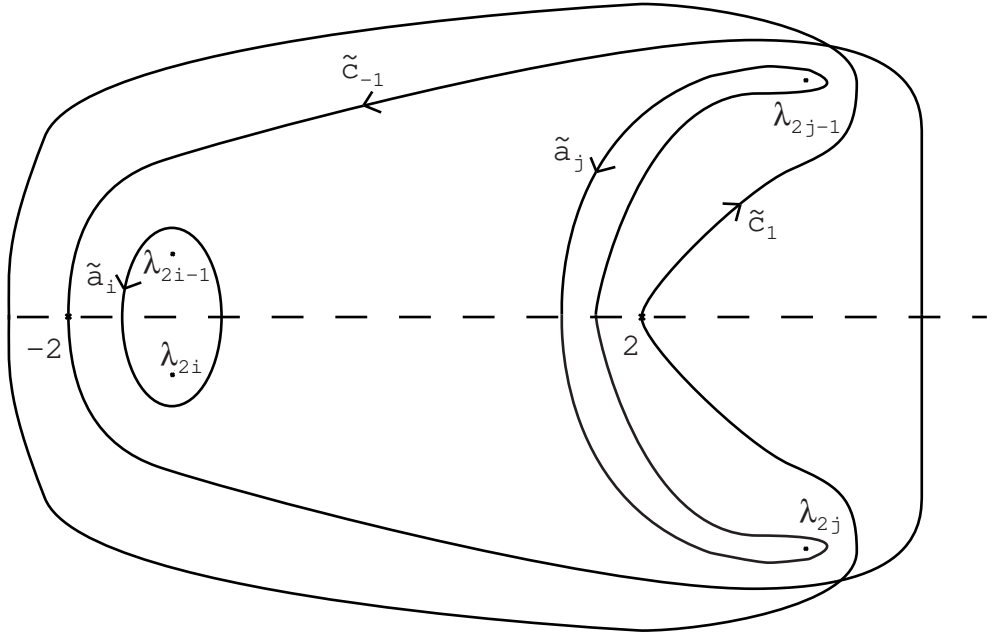}
\caption{The curves $\tilde a_i$ and $\tilde c_\pm$.\label{fig:evenCpm}}
\end{figure}

Choose lifts of the curves \mbox{$\tilde a_1,\ldots,\tilde a_n$} to $\Cpm$ and denote the homology classes of these lifts by \mbox{$a^\pm_1,\ldots,a^\pm_n$}. Choose also an open curve $c_1$ in $C_+$ covering $\tilde c_1$ and an open curve $c_{-1}$ in $C_-$ covering $\tilde c_{-1}$. We write $M_n$ for the space of $2n$-tuples $(\lambda_1,\ldots,\lambda_{2n})$ as above together with the choices described, and let $M_{n,\R}$ denote the subset of $M_n$ such that (see figure~\ref{fig:evenCpm}):
\begin{enumerate}

\item $\lambda_{2i}=\bar\lambda_{2i-1}$ for $i=1,\ldots,n$,

\item $\tilde a_i$ intersects the real axis exactly twice, both times in the interval \mbox{$(-2,2)$},

\item the lifts of  \mbox{$\tilde a_1,\ldots,\tilde a_n$} to $C_+$  
are chosen so that the point where $\tilde a_i$ intersects the $z$-axis with 
positive orientation is lifted to a point where $w_+$ is positive,

\item the lifts of  \mbox{$\tilde 
a_1,\ldots,\tilde a_n$} to $C_{-}$  are chosen so 
that the point where $\tilde a_i$ intersects the $z$-axis with positive orientation is lifted to a point where $\frac{w_-}{i}$ is positive,

\item $c_1$ begins at a point where $w_+$ is positive.

\item $c_{-1}$ begins at a point where $\frac{w_-}{i}$ is positive. 

\end{enumerate}

For $p\in M_{n,\R}$, we may, as before, construct \mbox{$\pi:X\ra\mathbb{C}P^1$} as the fibre product of \mbox{$\pi_{+}:C_+\ra\mathbb{C}P^1$} and 
\mbox{$\pi_{-}:C_-\ra\mathbb{C}P^1$}. $X$  is given by the equation 
$$y^2=x\prod_{i=1}^{n}(x-\alpha_i)(x-{\alpha_i}^{-1})(x-\bar{\alpha_i})(x-{\bar{\alpha_i}}^{-1})\mbox{, where }\alpha_i+{\alpha_i}^{-1}=\lambda_i,$$
and $\pi$ by
$$\pi:(x,y)\mapsto x.$$ These identifications occur via the maps  
$$q_\pm(x,y)=\left(x+\frac{1}{x},\frac{(x\pm 1)y}{x^{n+1}}\right)=(z,w_\pm).$$ 
$X$ has genus $2n$ and possesses the holomorphic involutions
$$\begin{array}{rccc}
i_{\pm}:&X&\longrightarrow&X\\
&(x,y)&\longmapsto&\left(\frac{1}{x},\frac{\pm y}{x^{2n+1}}\right).
\end{array}$$
The curves $\Cpm$ are the quotients of $X$ by these involutions, with 
quotient maps $q_\pm:X\ra\Cpm$, and the real structures $\rho_\pm$ on $\Cpm$ induce upon $X$ the real structure $$\rho:(x,y)\mapsto (\frac{1}{\bar x},\frac{\bar y}{\bar x^{2n+1}}).$$

\begin{figure}[h]
\hspace*{1.5cm}\includegraphics{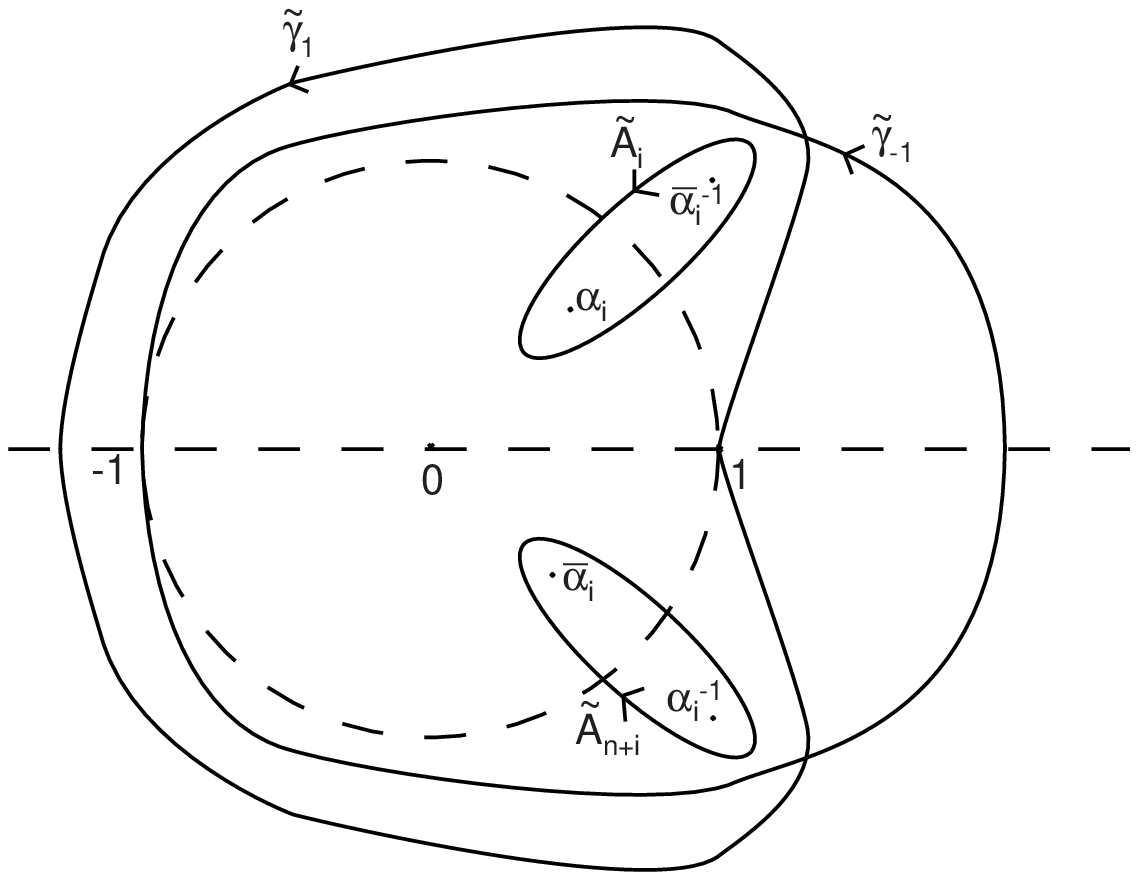}
\caption{The curves $\tilde A_i$ and $\tilde \gamma_\pm$.\label{fig:evenX}}
\end{figure}

For each $p\in M_{n,\R}$ there is a unique canonical basis 
\mbox{$A_1,\ldots,A_{2n}$}, \mbox{$B_1,\ldots,B_{2n}$} for the 
homology of $X$ such that \mbox{$A_1,\ldots,A_{2n}$} cover the 
homotopy classes of loops \mbox{$\tilde A_1,\ldots,\tilde A_{2n}$} 
shown in Figure~\ref{fig:evenX} and 
$$(q_\pm)_*(A_i)=\pm(q_\pm)_*(A_{n+i})=a^\pm_i.$$
There are also unique curves $\gamma_1$ and $\gamma_{-1}$ on $X$ such 
that $(q_\pm)_*(\gamma_{\pm 1})=c_{\pm 1}$; they project to 
$\tilde{\gamma_1}$ and $\tilde{\gamma_{-1}}$ of Figure~\ref{fig:evenX}.
For each $p\in M_n$ define differentials $\Wpm=\Wpm(p)$ on $\Cpm(p)$ by:
\begin{enumerate}
\item $\Wpm(p)$ are meromorphic differentials of the second kind: their 
only singularities are double poles at $z=\infty$, and they have no 
residues.
\item $\int_{a^\pm_i}\Wpm(p) = 0$ for 
$i=1,\ldots,n$.
\item As $z\ra\infty$, $\Wpm(p)\ra\frac{z^n dz}{w_\pm(p)}$. 
\end{enumerate}
In view of these defining conditions, we may write
$$\Wpm=\frac{\prod_{j=1}^n(z-\zeta^\pm_j)dz}{w_\pm}.$$
Define
$$I_{+}(p):=\left(\int_{c_1}\Omega_+(p),\int_{b^+_1}\Omega_+(p),\ldots,\int_{b^+_n}\Omega_+(p)\right),$$
$$I_{-}(p):=\im\left(\int_{c_{-1}}\Omega_-(p),\int_{b^-_1}\Omega_-(p),\ldots,\int_{b^-_n}\Omega_-(p)\right).$$

For $p\in M_{n,\R}$, $I_+(p)$ and $I_-(p)$ are real, since, writing ${\cal A}^\pm$ for the subgroups of \mbox{$H_1(\Cpm,\mathbb{Z})$} generated by the $a^\pm_i$, we then have
$$(\rho_\pm)_*(b^\pm_i)=b^\pm_i 
\mbox{ mod ${\cal A}^\pm$},\:  (\rho_\pm)_*(c_{\pm 1})=c_{\pm 1} \mbox{ 
mod ${\cal A}^\pm$},$$
and
$$\rho^*_\pm(\Wpm)=\pm\overline\Wpm.$$

\begin{thm}
    For each positive integer $m$ and integer $n$ with $0\leq n\leq m$ there exists $p\in M_{n,\R}$ such that 
\begin{enumerate}    
\item{$\zeta^+_j(p)$, $j=1,\ldots,n$ are pairwise distinct, as are 
$\zeta^-_j(p)$, $j=1\ldots n$}
\item $\mathbb R^{2n+2}$ is spanned by the vectors $(I_{+}(p),0)$, $(0,I_{-}(p))$ and $\pli\left(I_{+}(p), I_{-}(p)\right)$,\\ $i=1,\ldots,2n$.
\item $5\left(D_{+}(p)-D_{-}(p)\right)+4T^{k}_p\left(\eta^{+}(p)-\eta^{-}(p)\right)\neq 0$, for $0\leq k\leq m-n$,\\ where $D_\pm(p)$, $\eta^\pm(p)$ and $T_p$ are defined as follows:\\
Assume that \setcounter{enumi}{1}(\it\arabic{enumi}) and \setcounter{enumi}{2}(\it\arabic{enumi}) of this theorem hold for 
$p\in M_{n,\R}$. Then the differentials 
$\frac{\Wpm(p)}{z-\zeta^\pm_j}$ are a basis for the holomorphic 
differentials on $\Cpm(p)$. Thus we may define 
$c^\pm_j(p)$ by the equations
$$\frac{3}{2}\int_{a^\pm_i}z\Wpm(p)+\sum_{j=1}^{n}c^\pm_j(p)\int_{a^\pm_i}\frac{\Wpm(p)}{z-\zeta^\pm_j}=0,\:i,j=1,\ldots,n,$$
and set 
$$\Wpmh(p):=\frac{3}{2}z\Wpm(p)+\sum_{j=1}^{n}c^\pm_j(p)\frac{\Wpm(p)}{z-\zeta^\pm_j},$$ 
 $$\widehat{I}_+(p):=\left(\int_{c_1}\widehat{\Omega}_+(p),\int_{b^+_1}\widehat{\Omega}_+(p),\ldots,\int_{b^+_n}\widehat{\Omega}_+(p)\right),$$ 
$$\widehat{I}_-(p):=\im\left(\int_{c_{-1}}\widehat{\Omega}_-(p),\int_{b^-_1}\widehat{\Omega}_-(p),\ldots,\int_{b^-_n}\widehat{\Omega}_-(p)\right).$$ 
Then we define \mbox{$\eta^{\pm}(p)$} by\\
${\ds(\widehat{I}_+(p),\widehat I_-(p))- \eta^+(p)(I_+(p),0)-\eta^-(p)(0,I_-(p))}$
\vspace{-5mm}\begin{flushright}${\ds\in\mbox{span}\left\{\pli(I_+(p),I_-(p))\mbox{, }i=1\ldots n\right\},}$\end{flushright}
put 
$$D_\pm(p):=\frac{1}{2}\left(\mp 2 + \sum_{i=1}^{2n}\lambda_{i}\right)-\sum_{j=1}^n\zeta^\pm_j,$$ 
and let $T_p$ be the linear fractional transformation
$$T_p:x\mapsto (D_+(p)-D_-(p))\frac{3x-4(D_+(p)-D_-(p))}{4x-5(D_+(p)-D_-(p))}.$$  
\end{enumerate}
\label{thm}
\end{thm}

\noindent{\bf Proof:} For each fixed $m$, this theorem can be proven by induction. The induction step is both similar to and simpler than that detailed in the odd genus case, and is almost identical to the induction step described in \cite{EKT:93}. We thus omit it here, but demonstrate the existence of $p\in M_{0,\R}$ satisfying (1) and (2) and such that, for all $k\geq 0$, $5\left(D_{+}(p)-D_{-}(p)\right)+4T^{k}_p\left(\eta^{+}(p)-\eta^{-}(p)\right)\neq 0$. 

Consider the curves $C_\pm$ given by $$w_\pm^2=(z\pm 2).$$ Let $c_1$ be an open curve on $C_+$ from $(2,2)$ to $(2,-2)$, and $c_{-1}$ one on $C_-$ from $(-2,2\im)$ to $(-2,-2\im)$. We have 
$$\Wpm=\frac{dz}{w_\pm},$$
which gives $$I_\pm=\mp 8.$$
Also, $$\widehat\Omega_\pm=\frac{3}{2}z\Wpm,$$
which gives rise to $$\widehat I_\pm=8.$$
$\eta^\pm$ are defined by the equation 
$$(\widehat I_+,\widehat I_-)=\eta^+(I_+,0)+\eta^-(0,I_-),$$
so $$\eta^\pm=\mp 1.$$
We have $$D_\pm = \mp 1,$$ and hence the linear transformation $T$ is given by 
$$T:x\mapsto -\frac{3x+8}{2x+5}.$$
The statement $$5\left(D_{+}-D_{-}\right)+4T^{k}\left(\eta^{+}-\eta^{-}\right)=0$$ may be written as
$$T^k(-2)=\frac54,$$ but $-2$ is a fixed point of $T$, so this is false for all $k\geq 0$. \begin{flushright}$\Box$ \end{flushright}


\bibliographystyle{plain}
\bibliography{main_nocomments}
\end{document}